\documentclass[10pt,preprint]{elsarticle}

\usepackage{amssymb,amsfonts,amsmath,mathrsfs,mathtools}
\usepackage{graphicx,epsfig,subfigure}
\usepackage{bm}
\usepackage[margin=1in,papersize={8.5in,11in}]{geometry}
\usepackage{times}
\usepackage{multirow}
\usepackage{color}
\usepackage[title]{appendix}%
\usepackage{algorithm}%
\usepackage{algorithmicx}%
\usepackage{algpseudocode}%
\usepackage{listings}%

\def\vs{\vspace{0.2cm}}

\renewcommand{\top}{\text{T}}

\def\I{\mathcal{I}}

\def\J{\mathcal{J}}
\def\M{\mathcal{M}}

\def\r{\bm r}
\def\DEIM{\texttt{DEIM}}

\newtheorem{lemma}{\bf Lemma}[section]
\newtheorem{corollary}{\bf Corollary}[section]
\newtheorem{theorem}{\bf Theorem}[section]

\newtheorem{proposition}{\bf Proposition}[section]

\newenvironment{proof}{{\noindent \bf \em Proof:}}{\hfill$\square$}

\title{Collocation methods for nonlinear differential equations on low-rank manifolds}

\begin{document}
\begin{frontmatter}

\author[lbnl]{Alec Dektor\corref{correspondingAuthor}}
\ead{adektor@lbl.gov}

\address[lbnl]{Applied Mathematics \& Computational Research Division, Lawrence Berkeley National Laboratory, Berkeley (CA) 94720, USA.}

\cortext[correspondingAuthor]{Corresponding author}

\journal{Arxiv}

\begin{abstract}

We introduce new methods for integrating nonlinear differential equations on low-rank manifolds. 
These methods rely on interpolatory projections onto the tangent space, enabling low-rank time integration of vector fields that can be evaluated entry-wise. 
A key advantage of our approach is that it does not require the vector field to exhibit low-rank structure, thereby overcoming significant limitations of traditional dynamical low-rank methods based on orthogonal projection. 
To construct the interpolatory projectors, we develop a sparse tensor sampling algorithm based on the discrete empirical interpolation method (DEIM) that parameterizes tensor train manifolds and their tangent spaces with cross interpolation. 
Using these projectors, we propose two time integration schemes on low-rank tensor train manifolds. 
The first scheme integrates the solution at selected interpolation indices and constructs the solution with cross interpolation. The second scheme generalizes the well-known orthogonal projector-splitting integrator to interpolatory projectors. 
We demonstrate the proposed methods with applications to several tensor differential equations arising from the discretization of partial differential equations. 

\end{abstract}

\begin{keyword} 
low-rank approximation \sep
time-dependent tensors \sep
tensor differential equations \sep
tensor cross approximation \sep
tensor train format
\end{keyword}

\end{frontmatter}

\section{Introduction}
\label{sec:intro}

Consider the initial value problem 
\begin{equation}
\label{nonlinear_PDE0}
\displaystyle \frac{\partial u(\bm x,t)}{\partial t} = \mathcal{G}({u},\bm x, t),\qquad 
{u}(\bm x,0) = {u}_0(\bm x),
\end{equation}
governing the time evolution of a quantity of interest $u:\Omega \times [0,T] \to \mathbb{R}$, where $\Omega$ is a subset of $\mathbb{R}^d$ ($d \gg 1$) and $\mathcal{G}$ is a nonlinear operator that may depend on $\bm x$ and $t$. 
Equations of the form \eqref{nonlinear_PDE0} are found in many areas of physical sciences, engineering, and mathematics. 
For example, in applications of kinetic theory such as the Fokker--Planck equation \cite{Risken} and the Boltzmann equation \cite{cercignani1988}, in optimal mass transport \cite{Osher2019}, and as finite-dimensional approximations of functional differential equations \cite{Venturi2021,Venturi2018}. 
Discretizing \eqref{nonlinear_PDE0} with a method of lines yields the tensor differential equation 
%
%
%
\begin{equation}
\label{nonlinear_ODE0}
\displaystyle \frac{d X(t)}{dt} = G(X(t),t), \qquad 
X(0) = X_0,
\end{equation}
where $X(t):[0,T] \to \mathbb{R}^{n_1 \times \cdots \times n_d}$ is the time-dependent solution tensor and $G: \mathbb{R}^{n_1 \times \cdots \times n_d} \times [0,T] \to  \mathbb{R}^{n_1 \times \cdots \times n_d}$ is a discrete form of the operator $\mathcal{G}$. 
At any time $t$, the solution tensor $X(t)$ has $\mathcal{O}(n^d)$ degrees of freedom that make its computation and storage prohibitively expensive, even for small $d$. 

Several algorithms based on tensor networks have recently been proposed to reduce the number of degrees of freedom in the solution tensor $X(t)$ and integrate \eqref{nonlinear_ODE0} at a reasonable computational cost. 
A tensor network is a factorization of a high-dimensional tensor, such as $X(t)$, into a network of low-dimensional tensors with significantly fewer degrees of freedom. 
The number of degrees of freedom in a network depends on the chosen tensor format, e.g., 
tensor train (TT) \cite{OseledetsTT}, Tucker \cite{malik2018low,dolgov2021functional}, Hierarchical Tucker \cite{grasedyck2010hierarchical,lubich2013dynamical,grasedyck2018distributed} or canonical polyadic (CP) \cite{kolda2009tensor}, 
and the tensor rank. 
For instance, tensors in the TT format with rank $r$ can be parameterized with $\mathcal{O}(dnr^2)$ degrees of freedom, a significant reduction from $\mathcal{O}(n^d)$ when the rank $r$ is sufficiently small. 
The set of all tensors in a chosen format with fixed rank forms a smooth manifold on which the solution to \eqref{nonlinear_ODE0} can be integrated. 

Two classes of algorithms for integrating \eqref{nonlinear_ODE0} on smooth tensor manifolds are step-truncation methods \cite{rodgers2020step-truncation} and dynamical low-rank methods \cite{Dektor_dyn_approx}. 
Step-truncation methods allow the solution rank to naturally increase in a controlled manner during a time step before truncating back to the desired rank. 
Dynamical low-rank methods integrate the solution on a fixed-rank manifold by projecting $G(X,t)$ onto a tangent space of the manifold at each time $t$. 
Both methods aim to efficiently compute the best approximate solution to \eqref{nonlinear_ODE0} on a fixed-rank tensor manifold at each time $t$ and are consistent with each other as the temporal step-size approaches zero (see, e.g., \cite[Section 3.3]{adaptive_rank}). 
These methods have been utilized for several applications including uncertainty quantification \cite{babaee2017robust,sapsis2009dynamically}, plasma physics \cite{ye2023quantized,einkemmer2018low}, numerical approximation of functional differential equations \cite{Venturi2021,rodgers2024tensor} and machine learning \cite{schotthofer2022low,savostianova2024robust}, 
and substantial research efforts have recently been made to improve their accuracy, efficiency, and robustness. 
These efforts have resulted in several innovations including rank-adaptive integrators \cite{adaptive_rank,ceruti2022rank,rodgers2020step-truncation}, implicit low-rank methods \cite{Rodgers_implicit,Nakao2023,Sutti2023}, conservative low-rank methods \cite{guo2024conservative,baumann2024energy,einkemmer2023robust,einkemmer2023conservation} and coordinate-adaptive low-rank methods \cite{Coordinate_flows,Coordinate_adaptive_integration}.

Despite these recent advancements, existing low-rank time integration schemes still have significant limitations in their applicability. 
Most notably, they require the tensor-valued map $G$ defining the differential equation \eqref{nonlinear_ODE0} to have low-rank structure complementary to that of the solution. 
Such low-rank structure is key to increasing the solution rank in a controlled manner for step-truncation methods or efficiently computing the orthogonal projection of $G(X,t)$ onto the tangent space for dynamical low-rank methods. 
In either case, the low-rank structure of $G$ is crucial for obtaining practical time integration schemes with computational cost and storage requirements comparable to the storage cost of the chosen tensor format. 
However, many instances of $G$ that arise from discretizing \eqref{nonlinear_PDE0} lack low-rank structure. 
For example, when $G$ includes a polynomial nonlinearity computing a low-rank representation of $G(X,t)$ or its orthogonal projection onto the tangent space is expensive due to the non-optimal rank that results from multiplying low-rank tensors. 
The situation is worse for other common nonlinearities, such as exponential or fractional, as there are currently no reliable algorithms for performing these nonlinear arithmetic operations with low-rank tensors. 
In such cases, \eqref{nonlinear_ODE0} may admit an approximate low-rank solution. However, existing step-truncation and dynamical low-rank methods cannot efficiently compute it.

In this paper, we introduce a new class of dynamical low-rank methods that can efficiently integrate the solution to \eqref{nonlinear_ODE0} on a low-rank tensor manifold even when $G$ does not have low-rank structure. 
Our proposed methods rely on a new class of oblique projectors onto low-rank tangent spaces with a cross interpolation property. 
In the context of dynamical low-rank approximation, these projectors collocate \eqref{nonlinear_ODE0} on a tensor manifold 
and yield equations of motion on the manifold that are efficient to integrate whenever it is possible to evaluate $G$ entry-wise. 
The oblique projectors are defined by sets of multi-indices that identify tensor fibers along which the projector interpolates. 
To select these multi-indices, we introduce a new algorithm, based on the discrete empirical interpolation method (DEIM), to efficiently compute indices that parameterize tensor manifolds with cross interpolation \cite{TT-cross}. 
Using these projectors we propose two low-rank time integration schemes. 
The first integrates subtensors of the solution defining a tensor cross interpolant and then constructs the low-rank solution later in time with tensor cross interpolation. 
The second we obtain by applying a splitting scheme to the oblique tangent space projector. 
Splitting schemes for orthogonal tangent space projectors were introduced in \cite{lubich2014projector} for matrices and were subsequently generalized to TTs \cite{Lubich_2015}, Tucker tensors \cite{LubichTucker} and tree tensor networks \cite{Lubich_tree_tensor}. 
Our method directly generalizes these orthogonal projector-splitting schemes to oblique projectors in the TT format. 
%
%
Related time integration schemes based on oblique projections onto the low-rank manifold (instead of its tangent space) were recently proposed for computing low-rank approximations to matrix differential equations \cite{Oblique-CUR, Naderi2023} and tensor differential equations \cite{ghahremani2024deim,ghahremani2024cross}.

The rest of this paper is organized as follows. 
In Section \ref{sec:matrix_DLR} we introduce interpolatory tangent space projectors and the proposed dynamical low-rank methods on matrix manifolds ($d=2$). 
In Section \ref{sec:TT_mfld} we briefly recall the TT format and orthogonalization of tensors in the TT format. 
In Section \ref{sec:projections} we recall the orthogonal projector onto the TT tangent space and introduce new oblique projectors onto the tangent space. 
Then we describe a special class of oblique projectors with a cross interpolation property. 
In Section \ref{sec:index_selection} we present a new index selection algorithm, referred to as TT-cross-DEIM, for constructing oblique projectors onto the TT tangent space. 
We show that indices obtained with the TT-cross-DEIM algorithm define oblique tangent space projectors and parameterize TT manifolds with cross interpolation. 
In Section \ref{sec:tt_integration} we introduce new dynamical low-rank time integration schemes for \eqref{nonlinear_ODE0} using oblique tangent space projectors. 
In Section \ref{sec:numerics} we demonstrate the proposed dynamical low-rank methods and compare the results with existing time integration methods on low-rank tensor manifolds. 
The main findings are summarized in Section \ref{sec:conclusions}. 

\section{Dynamical low-rank matrix approximation} 
\label{sec:matrix_DLR}

\begin{figure}[!t]
\centering
\includegraphics[scale=0.6]{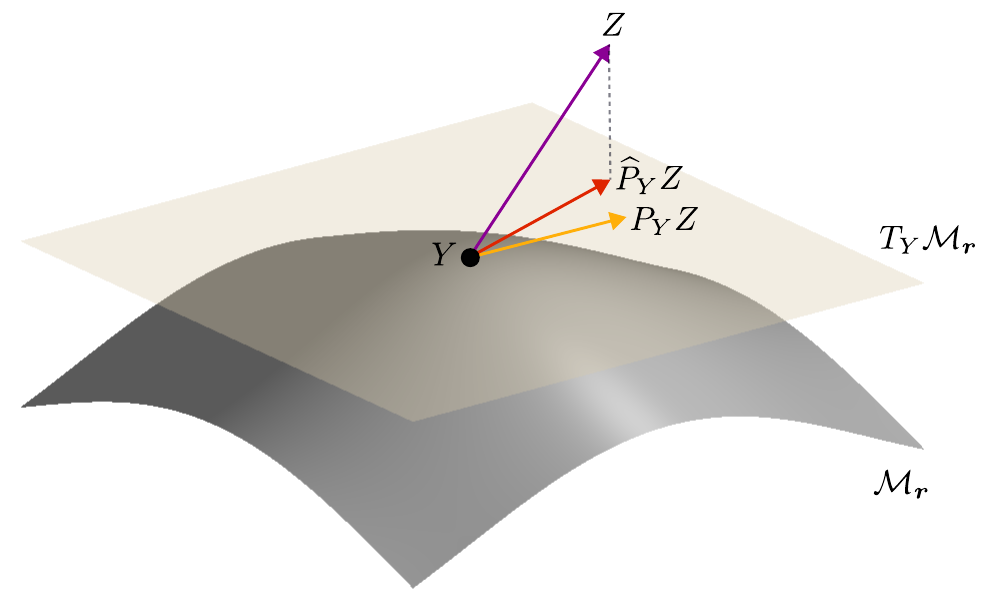}
\caption{
A sketch of the low-rank manifold $\M_{\r}$ and its tangent space $T_{Y}\M_{\r}$ at the point $Y \in \M_{\r}$. 
Also depicted is $Z \in \mathbb{R}^{n_1 \times \cdots \times n_d}$ and its orthogonal projection $\widehat{P}_{Y} Z$ and oblique projection $P_{Y} Z$ onto the tangent space $T_{Y}\M_{\r}$. 
The orthogonal projection is the best approximation of $Z$ on the tangent space with respect to the Frobenius norm but is impractical compute when $G$ lacks rank structure. 
The oblique projection is a quasi-optimal approximation of $Z$ on the tangent space that is efficient to compute for any $G$ that can be evaluated entry-wise. 
Such oblique projectors allow us to efficiently apply dynamical low-rank methods to a broad class of nonlinear differential equations. 
}
\label{fig:tensor_mfld}
\end{figure}

Before introducing dynamical low-rank tensor approximation, we describe the proposed low-rank methods for matrices ($d=2$). The goal is to find an approximate solution $Y(t)$ to \eqref{nonlinear_ODE0} that lies on the manifold $\M_r$ of rank-$r$ matrices for all time $t$. A rank-$r$ matrix can be expressed through left and right factor matrices with dimensions $n_1r$ and $n_2r$, respectively. Thus approximating the solution to \eqref{nonlinear_ODE0} with $Y(t) \in \M_r$ reduces the storage cost from $n_1 n_2$ to $(n_1 + n_2)r$. Dynamical low-rank methods integrate the approximate solution $Y(t)$ on the manifold $\M_r$ by projecting \eqref{nonlinear_ODE0} onto a tangent space of $\M_r$ at each time $t$. Assuming the initial condition $X_0$ belongs to $\M_r$, this procedure yields the evolution equation 
\begin{equation}
\label{DLRA}
\frac{dY(t)}{dt} = P_{Y(t)} G(Y(t),t), \qquad Y(0) = X_0, 
\end{equation}
where $P_{Y(t)}:\mathbb{R}^{n_1\times n_2} \to T_{Y(t)}\M_r$ projects onto the tangent space of $\M_r$ at $Y(t)$. The solution to \eqref{DLRA} remains on $\M_r$ for all $t \geq 0$ and serves as an approximate solution to \eqref{nonlinear_ODE0}. Classical dynamical low-rank methods use an orthogonal tangent space projector to minimize the error of the approximation in the Frobenius norm at each time $t$. Let $U(t)$ and $V(t)$ be matrices whose columns form orthonormal bases for the range and co-range of $Y(t)$, respectively. Such matrices can be obtained, for example, from the SVD of $Y(t)$. Then the orthogonal tangent space projector can be expressed as \cite{lubich2014projector} 
\begin{equation}
\label{orth_TS_proj}
\widehat{P}_Y Z = Z VV^{\top}
- UU^{\top} Z VV^{\top} 
+ UU^{\top} Z, 
\end{equation}
where $Z \in \mathbb{R}^{n_1 \times n_2}$ and we suppressed dependence on $t$ for simplicity. When $Z = G(Y(t),t)$ computing the projection \eqref{oblique_TS_proj} at each time step to integrate \eqref{DLRA} can be computationally expensive, especially if $G$ does not have low-rank structure. For most nonlinear functions $G$, computing such orthogonal projection has computational cost scaling as $n_1n_2$, making the integration of the approximate solution $Y(t) \in \M_r$ as expensive as solving for the full solution $X(t)$ using standard methods. 

\subsection{Interpolatory dynamical low-rank approximation}

\begin{algorithm}
 \caption{DEIM index selection (adapted from \cite{DEIM-CUR})}
\label{alg:DEIM}
\begin{algorithmic}[1]
\Require 
    $\bm V \in \mathbb{R}^{n \times r}$ with $n \geq r$
\Ensure 
	$\bm l$, a vector with $r$ distinct indices from $\{1,\ldots,n\}$
%
%
\State $\bm v = \bm V(:,1)$
\State $[~,l_1] = \max(|\bm v|)$
\For{$j = 2,3,\ldots,r$ }
        \State $\bm v = \bm V(:,j)$ 
        \State $\bm c = \bm V(l,1:j-1)^{-1} \bm v(\bm l)$
        \State $\bm r =\bm v - \bm V(:,1:j-1) \bm c$ 
        \State $[~,l_j] = \max(|\bm r|)$
        \State $\bm l = [\bm l; l_j]$
\EndFor
\end{algorithmic}
\end{algorithm}

To develop efficient dynamical low-rank integrators for nonlinear $G$ we propose a new interpolatory tangent space projector $P_Y$ to replace the orthogonal projector. 
This interpolatory projector is a specialized oblique tangent space projector, obtained by replacing the orthogonal projectors $UU^{\top}$ and $VV^{\top}$ in \eqref{orth_TS_proj} with oblique projectors onto the same spaces. A general form of these oblique projectors are 
$
U(A^{\top}U)^{-1}A^{\top}$ and $V(B^{\top}V)^{-1}B^{\top}$  
where $A\in \mathbb{R}^{n_1 \times r}$ and $B\in \mathbb{R}^{n_2 \times r}$ are any matrices such that $(A^{\top}U)$ and $(B^{\top}V)$ are invertible. 
Interpolatory projectors onto the columns of $U$ and $V$ are obtained by selecting $A$ and $B$ as specific columns of the identity matrix with appropriate dimensions. Specifically $A = I_{n_1}(:,\I)$ and $B = I_{n_2}(:,\J)$ where $\I$ contains $r$ indices from $\{1,\ldots,n_1\}$ and $\J$ contains $r$ indices from $\{1,\ldots,n_2\}$. The matrix $A$ can extract $r$ rows from $Z$ via left multiplication $A^{\top}Z = Z(\I,:)$ and the matrix $B$ can extract $r$ columns from $Z$ via right multiplication $ZB = Z(:,\J)$. Moreover, $A$ can be defined by choosing $r$ indices in the set $\I$, and similarly, $B$ can be defined by choosing $r$ indices in the set $\J$. Such indices are to be chosen so that $A^{\top}U = U(\I,:)$ and $B^{\top} V = V(\J,:)$ are invertible, and ideally with a small condition number.  
This can be achieved with a sparse sampling algorithm such as the discrete empirical interpolation method (DEIM). The DEIM, summarized in Algorithm~\ref{alg:DEIM}, is a greedy algorithm that selects an index for each column of $U$ or $V$ to minimize the condition number of the interpolatory projector as much as possible. Other sparse sampling strategies, such as Q-DEIM or oversampling methods, can also be used and may yield better-conditioned interpolatory projectors than DEIM in certain cases. 
Replacing the orthogonal projectors in \eqref{orth_TS_proj} with interpolatory projectors onto the columns of $U$ and $V$ results in an interpolatory projector onto the tangent space $T_{Y}\M_r$ of the form 
\begin{equation}
\label{oblique_TS_proj} 
P_Y Z = Z(:,\J)V(\J,:)^{-\top}V^{\top} 
- U U(\I,:)^{-1} Z(\I,\J) V(\J,:)^{-\top}V^{\top} 
+ U U(\I,:)^{-1} Z(\I,:). 
\end{equation}
It is easy to verify that $(P_YZ)(i,j) = Z(i,j)$ whenever $i \in \I$ or $j \in \J$. In other words, the projection $P_YZ$ interpolates $Z$ along rows and columns specified by the index sets $\I$ and $\J$.

\subsubsection{Matrix cross integrator} 
\label{sec:matrix_cross}
A straightforward low-rank time integration scheme can be derived by evaluating the dynamical low-rank evolution equation \eqref{DLRA} at the sampled indices $\I$ and $\J$ and leveraging the interpolation property. This leads to the system of evolution equations
\begin{equation}
\label{matrix_cross_ev}
\begin{aligned}
    \frac{d Y(\I(t),:,t)}{dt} &= G_Y(\I(t),:,t), \\
    \frac{d Y(:,\J(t),t)}{dt} &= G_Y(:,\J(t),t),
\end{aligned}
\end{equation}
where $G_Y(t)=G(Y(t),t)$. Equation \eqref{matrix_cross_ev} consists of $(n_1 + n_2)r$ coupled nonlinear differential equations governing the evolution of a subset of the entries in the approximate solution, which can be integrated using standard explicit or implicit methods. The indices $\I(t)$ and $\J(t)$ that define the interpolatory projector \eqref{oblique_TS_proj} are chosen at each time $t$ to ensure that the projector remains well-defined during time integration. 
If $G$ arises from the spatial discretization of a PDE \eqref{nonlinear_PDE0} involving differential operators, evaluating $G_Y\left(\I,:,t\right)$ and $G_Y\left(:,\J,t\right)$ requires entries of $Y$ at indices adjacent to $\I$ and $\J$. 
%
%
The values of $Y$ at adjacent indices can always be obtained by constructing the solution $Y(t)$ as a low-rank matrix using CUR decomposition. 
For example given $Y(\I,:)$ and $Y(:,\J)$ at any time $t$ the approximate solution $Y$ can be obtained as 
\begin{equation} \label{CUR1}
    Y = Y(:,\J) Y(\I,\J)^{-1} Y(\I,:).
\end{equation}
If $Y = U\Sigma V^{\top}$ is the SVD and the indices $\I$ and $\J$ defining the interpolatory projector are chosen so that $U(\I,:)$ and $V(\J,:)$ are well-conditioned, then constructing $Y$ using \eqref{CUR1} involves $Y(\I,\J)^{-1} = V(\J,:)^{-\top}\Sigma^{-1}U(\I,:)^{-1}$. Thus the condition number of the middle matrix is inversely proportional to the smallest singular value of $Y$. 
Instead of constructing $Y$ in this way, if we first take a QR decomposition $Y(:,\J) = QR$, then we can write $Y(\I,\J) = Q(\I,:)R$ and the CUR formula \eqref{CUR1} becomes
\begin{equation} \label{CUR2}
    Y = QQ(\I,:)^{-1}Y(\I,:), 
\end{equation}
which is an interpolatory projection of $Y(\I,:)$ onto the orthonormal basis $Q$. In particular the stability of constructing the solution using \eqref{CUR2} depends on the condition of the interpolatory projector and is independent of the singular values of $Y$. 
Computing the right hand side of \eqref{matrix_cross_ev}, and hence integrating the system, is efficient for any nonlinear $G$ that can be evaluated entry-wise, regardless of its low-rank structure. This enables us to efficiently compute dynamical low-rank approximations for problems where orthogonal tangent space projection is too expensive.

\subsubsection{Projector-splitting integrator}
\label{sec:proj_split_mat}
An alternative approach to the matrix-cross integrator presented above for integrating \eqref{DLRA} is to apply a standard splitting method, similar to the integrators proposed for orthogonal projectors in \cite{lubich2014projector}. Since \eqref{DLRA} is a sum of three terms, applying Lie-Trotter splitting yields three substeps commonly referred to as \textbf{K-}, \textbf{S-}, and \textbf{L-step}. 
Beginning from the rank-$r$ decomposition of the approximate solution $Y(t_0) = U(t_0) S(t_0) V^{\top}(t_0)$ at time $t_0$, each step updates a single factor matrix. After all three steps we obtain the low-rank factors for the approximate solution $Y(t_1) = U(t_1) S(t_1) V^{\top}(t_1)$ at time $t_1 = t_0 + \Delta t$. The substeps for interpolatory projector-splitting integrator are as follows. We denote by $\texttt{DEIM}(\cdot)$ a subroutine that takes a matrix of size $n \times r$ as input and outputs a collection of $r$ indices computed with the DEIM index selection (Algorithm~\ref{alg:DEIM}). 
\begin{itemize}
\item[1.] {\bf K-step}: update $U(t_0) \to U(t_1)$ and $S(t_0) \to R(t_1)$. 

\vs
Compute interpolation indices $\J = \DEIM(V(t_0))$. Then integrate the $n_1 \times r$ differential equation 
\begin{equation}
\label{K}
\frac{d  K(t)}{dt} = 
 G_{K}\left(:,\J,t\right) \left[ V\left(\J,:,t_0\right)\right]^{-\top}, \qquad  K(t_0) = U(t_0) S(t_0), 
\end{equation}
where $ G_{K}(t) = G\left( K(t)V(t_0)^{\top},t\right)$ from $t_0$ to $t_1$, and perform a QR-decomposition $ K(t_1) =  U(t_1)  R(t_1)$. 

\item[2.] {\bf S-step}: update $ R(t_1) \to \tilde{ S}(t_1)$. 

\vs
Compute interpolation indices $\I = \DEIM( U(t_1))$. 
Then integrate the $r \times r$ matrix differential equation 
\begin{equation}
\frac{d\tilde{ S}(t)}{dt} = -\left[ U\left(\I,:,t_1\right)\right]^{-1} 
 G_{ S}\left(\I,\J,t\right) \left[ V\left(\J,:,t_0\right)\right]^{-\top}, 
\qquad \tilde{ S}(t_0) =  R(t_1), 
\end{equation}
where $ G_{ S}(t) = G\left( U(t_1)  S(t)  V(t_0)^{\top},t\right)$ from $t_0$ to $t_1$. 

\item[3.] {\bf L-step}: update $ V(t_0) \to  V(t_1)$ and $\tilde{ S}(t_1) \to  S(t_1)$. 

\vs
Integrate the $n_2 \times r$ matrix differential equation 
\begin{equation}
\label{Lstep}
\frac{d  L(t)}{dt} = 
 G_{ L}\left(\I,:,t\right)^{\top} 
\left[ U\left(\I,:,t_1\right)\right]^{-\top}, \qquad  L(t_0) =  V(t_0)\tilde{ S}(t_1)^{\top}, 
\end{equation}
where $ G_{ L}(t) = G \left( U(t_1) L(t)^{\top},t\right)$ from time $t_0$ to $t_1$, 
and perform a QR-decomposition $ L(t_1) = V(t_1)  S(t_1)^{\top}$. 
\end{itemize}

\noindent
Similar to the matrix-cross integrator, the interpolatory projector-splitting integrator requires evaluating the output of $G$ at only a subset of $nr$ indices in the {\bf K-} and {\bf L-step} and $r^2$ indices in the {\bf S-step}. These evaluations can be performed efficiently for any $G$ that can be evaluated entry-wise. 
Meanwhile the corresponding steps of the orthogonal projector-splitting integrator (summarized in \cite{Ceruti2022}) involve inner products involving the output of $G$ that are only efficient to compute when $G$ has low-rank structure. 
The difference between the interpolatory projector-splitting integrator and the matrix cross integrator is the order in which interpolatory projection and time integration are performed. The former applies the interpolatory projection on the vector-field and then integrates the factor matrices of the solution. The latter integrates the solution at the specified indices and then performs an interpolatory projection onto the updated basis $Q$ in \eqref{CUR2}.

\section{Tensor train (TT) format}
\label{sec:TT_mfld} 

In the following sections we propose dynamical low-rank approximation with interpolatory projections for tensors $(d\geq 2)$ in the tensor train (TT) format. When $d=2$ the TT methods described hereafter reduce to the matrix methods discussed in Section~\ref{sec:matrix_DLR}. 
We begin with a brief review of the TT format and orthogonal TT representations. For a more detailed introduction to the TT format we refer the reader to \cite{OseledetsTT}. Throughout the remainder of this paper, matrices are denoted by boldface letters, while tensors are denoted by regular (non-bold) letters. 
The $k$th unfolding of a tensor $Y \in \mathbb{R}^{n_1 \times \cdots \times n_d}$ is the matrix 
$\bm Y^{\langle k \rangle} \in \mathbb{R}^{(n_1\cdots n_k)\times (n_{k+1}\cdots n_d)}$ 
with rows and columns indexed colexicographically. 
%
%
%
%
The TT-rank of $Y$ is defined as the vector $\bm r = (1,r_1,\ldots,r_d,1)$ where $r_k$ is the rank of the unfolding matrix $\bm Y^{\langle k \rangle}$. 
Any tensor $Y$ with TT-rank $\bm r$ can be represented in the TT format as 
\begin{equation}
\label{TT_SVD}
Y(i_1, i_2, \ldots, i_d) = C_1(i_1) C_2(i_2) \cdots C_d(i_d),
\end{equation}
where each $C_k$ is a $r_{k-1} \times n_k \times r_{k}$ tensor referred to as a TT-core and $C_k(i_k)$ is a $r_{k-1} \times r_k$ matrix for a fixed index $i_k$. 
Each TT-core has a left unfolding matrix $\bm{C}_k^{\langle l \rangle} \in \mathbb{R}^{r_{k-1} n_k \times r_k}$ and a right unfolding $\bm{C}_k^{\langle r \rangle} \in \mathbb{R}^{r_{k-1} \times n_k r_k}$ obtained by reshaping the elements of $C_k$ 
$$\bm C_k^{\langle l \rangle}(\alpha_{k-1}i_k,\alpha_k) = 
\bm C_k^{\langle r \rangle}(\alpha_{k-1}, i_k \alpha_k) = 
C_k(\alpha_{k-1},i_k,\alpha_k),
\qquad k =1,2,\ldots,d. $$ 
The left and right unfoldings of each TT-core is full rank whenever $Y$ has TT-rank $\bm r$. 
To simplify notation of tensors in the TT format we often omit indices so that \eqref{TT_SVD} is replaced by $Y = C_1 C_2 \cdots C_d$. 
To obtain even more compact representations, we define partial products of TT-cores $C_{\leq k} \in \mathbb{R}^{n_1 \times \cdots \times n_k \times r_k}$ and $C_{> k} \in \mathbb{R}^{r_k \times n_{k+1} \times \cdots \times n_d}$ with entries 
\begin{equation}
\label{partial_tensor}
\begin{split}
     C_{\leq k}(i_1,\ldots,i_k,:) &= C_1(i_1) \cdots C_k(i_k,:), 
    \\
     C_{>k}(:,i_{k+1},\ldots,i_d) &= C_{k+1}(:,i_{k+1}) \cdots C_d(i_d), 
\end{split}
\end{equation}
so that $Y = C_{\leq k} C_{>k}$. 
We also define certain unfolding matrices of these partial product tensors 
\begin{equation}
\label{partial_unfolding}
\begin{split}
    \bm C_{\leq k}(i_1\cdots i_k,:) &= C_{\leq k}(i_1,\ldots,i_k,:), 
    \\
    \bm C_{>k}(i_{k+1} \cdots i_d,:) &= C_{>k}(:,i_{k+1},\ldots,i_d),
\end{split}
\end{equation}
where $\bm C_{\leq k} \in\mathbb{R}^{(n_1\cdots n_k)\times r_k}$ and $\bm C_{> k} \in\mathbb{R}^{(n_{k+1}\cdots n_d)\times r_k}$, which allows us to write the $k$th unfolding matrix of $Y$ as $\bm Y^{\langle k \rangle} = \bm C_{\leq k} \bm C_{>k}^{\top}$. 

\begin{figure}[!t]
\centerline{\footnotesize\hspace{1cm} (a)   \hspace{6cm} (b) \hspace{1.6cm}}
\centering
\includegraphics[scale=0.5]{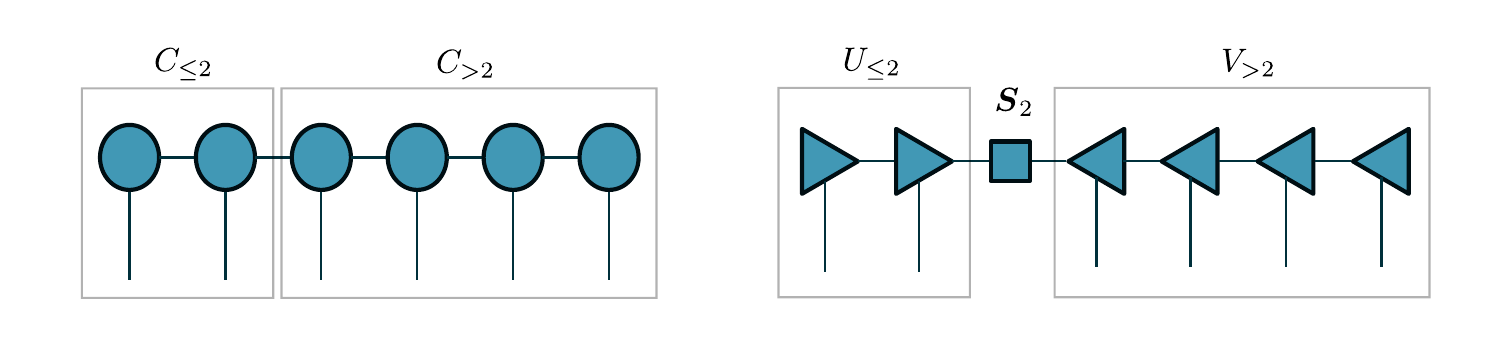}
\caption{
Tensor network diagrams of a $d=6$ dimensional tensor in the TT format. 
(a) No orthogonalization with the TT-core partial products $C_{\leq 2}$ and $C_{>2}$ indicated. 
(b) Orthogonalized TT \eqref{orth_TT} with $k=2$ and left orthogonal TT-cores $U_{\leq 2}$ and right orthogonal TT-cores $V_{>2}$ indicated. 
}
\label{fig:TT}
\end{figure}

\subsection{Orthogonalization of tensor trains}
\label{sec:orth_TT}

Orthogonal TT representations are fundamental for executing many operations in the TT format. 
We will use them in this paper to obtain projectors onto the tangent spaces of TT manifolds. 
Hereafter, we recall an algorithm for orthogonalizing TTs by recursively applying QR-decomposition to TT-core unfoldings. 
Begin by taking a QR-decomposition of the left unfolding of $ C_1$ 
%
$$
    \bm C_1^{\langle l \rangle} = \bm U_1^{\langle l \rangle} \bm R_1, 
$$
%
to obtain the matrix $\bm R_1 \in \mathbb{R}^{r_1 \times r_1}$ and the new TT-core $U_1 \in \mathbb{R}^{r_0 \times n_1 \times r_1}$ defined by its left unfolding. 
The new TT-core is called left-orthogonal because it satisfies 
$$
\left[\bm U_1^{\langle l \rangle}\right]^{\top} \bm U_1^{\langle l \rangle} = \bm I_{r_1}, 
$$
where $\bm I_{r_1}$ denotes the $r_1 \times r_1$ identity matrix. 
Next define a new second core $\widehat{C}_2(i_2) = \bm R_1 C_2(i_2)$ to obtain the TT representation 
$Y = U_1 \widehat{C}_2 C_{>2}$ 
of the tensor \eqref{TT_SVD}.
Then take a QR-decomposition of the left unfolding of the second core 
$$
    \widehat{\bm C}_2^{\langle l \rangle}
    = \bm U_2^{\langle l \rangle} \bm R_2,
$$
to obtain $\bm R_2 \in \mathbb{R}^{r_2 \times r_2}$ and the left-orthogonal TT-core $U_2 \in \mathbb{R}^{r_1 \times n_2 \times r_2}$. 
Define a new third core $\widehat{C}_3(i_3) = \bm R_2 C_3(i_3)$ to write $Y = U_{\leq 2}  \widehat{C}_3 C_{>3}$ where now the first two cores are left orthogonal. 
Proceeding recursively in this way we obtain the TT representation 
$$
    Y = U_{\leq k} \bm R_k C_{>k},
$$
where $U_j$ is left-orthogonal for $j = 1,2,\ldots,k$ and $\bm R_k \in \mathbb{R}^{r_k \times r_k}$. 
Similar to orthogonalizing cores from left to right, we can also orthogonalize cores from right to left by recursively performing QR-decompositions on right unfoldings (see \cite[Section 3]{OseledetsTT}) to obtain 
$$
Y = U_{\leq k} \bm R_k \bm R_{k+1} V_{>k},
$$
where $\bm R_{k+1} \in \mathbb{R}^{r_k \times r_k}$ and the $V_{j}$ are right-orthogonal, i.e., 
$$
    \bm V_j^{\langle r \rangle} \left[ \bm V_j^{\langle r \rangle}\right]^{\top} = \bm I_{r_{j-1}}, \qquad j = k+1,\ldots,d. 
$$
Letting $\bm S_k = \bm R_k \bm R_{k+1}$ we obtain the orthogonalized TT representation 
\begin{equation}
    \label{orth_TT}
Y = U_{\leq k} \bm S_k V_{>k}. 
\end{equation}
Utilizing the unfolding matrices of partial products \eqref{partial_unfolding} we also have a decomposition of the $k$th unfolding matrix 
\begin{equation}
\label{svd_k}
   \bm Y^{\langle k \rangle} =  \bm U_{\leq k} \bm S_k  \bm V_{>k}^{\top}. 
\end{equation}
It follows from the left orthogonality of $U_j$ and right orthogonality of $V_j$ that the columns of $\bm U_{\leq k}$ and $\bm V_{>k}$ are orthonormal, i.e., 
$\bm U_{\leq k}^{\top} \bm U_{\leq k} 
= \bm V_{>k}^{\top}\bm V_{>k} 
= \bm I_{r_k}$. 
The decomposition \eqref{svd_k} resembles a SVD however $\bm S_k$ is not necessarily diagonal. 
If the TT-rank of $Y$ is $\bm r$ then $\bm S_k$ is invertible. 
What is important for the projectors defined in the subsequent section is that the columns of $\bm U_{\leq k}$ and $\bm V_{>k}$ form orthonormal bases for the range and co-range of $\bm Y^{\langle k \rangle}$ respectively.

\section{Projections onto the TT tangent space}
\label{sec:projections}

It is well-known that the collection of all rank-$\bm r$ TTs 
\begin{equation}
\label{TT-mfld}
\mathcal{M}_{\bm r} = \{Y \in \mathbb{R}^{n_1\times \cdots \times n_d}  \mid  \text{TT-rank}(Y) = \bm r\},
\end{equation}
is a smooth embedded submanifold of $\mathbb{R}^{n_1 \times \cdots \times n_d}$ \cite{TT-mfld}. 
Hence for any tensor $Y \in \M_{\r}$ we can define the tangent space $T_{Y}\M_{\r}$, which is 
%
a vector subspace of $\mathbb{R}^{n_1 \times \cdots \times n_d}$ that linearizes the manifold $\M_{\r}$ around $Y$. 
Given a TT representation \eqref{TT_SVD} of $Y$, any element of the tangent space can be written (non-uniquely) as 
\begin{equation}
\label{tangent_element}
    \delta Y = \delta C_1 C_2 \cdots C_d + C_1 \delta C_2 C_3 \cdots C_d + \cdots + C_1 \cdots C_{d-1} \delta C_d, 
\end{equation}
where $\delta C_k \in \mathbb{R}^{r_{k-1} \times n_k \times r_k}$ are first order variations of the TT-cores. 

\subsection{The orthogonal tangent space projector} 
The orthogonal projector 
$\widehat{P}_{Y}:\mathbb{R}^{n_1 \times \cdots \times n_d} \to T_{Y} \M_{\r}$ onto the tangent space considered in \cite{Lubich_2015} determines the best approximation of a given tensor $Z \in \mathbb{R}^{n_1 \times \cdots \times n_d}$ in the tangent space relative to the Frobenius norm. Such orthogonal projector can be constructed from the orthogonal projectors 
\begin{equation}
\label{orthog_range_co}
\begin{split}
    \widehat{\bm P}_{\leq k} = \bm U_{\leq k} \bm U_{\leq k}^{\top}, \qquad 
    \widehat{\bm P}_{>k} = \bm V_{>k} \bm V_{> k}^{\top}, \qquad k = 1,2,\ldots,d-1, 
\end{split}
\end{equation}
onto the range and co-range of $\bm Y^{\langle k \rangle}$. The orthogonal bases $\bm U_{\leq k}$ and $\bm V_{>k}$ for the range and co-range of $\bm Y^{\langle k \rangle}$ can be obtained from the TT-orthogonalization procedure described in Section \ref{sec:orth_TT}. 
These projectors act on the matrix space $\mathbb{R}^{(n_1 \cdots n_k) \times (n_{k+1} \cdots n_d)}$. 
To construct the orthogonal tangent space projector $\widehat{P}_Y$ it is convenient to define projectors corresponding to \eqref{orthog_range_co} that act on the tensor space $\mathbb{R}^{n_1  \times \cdots \times n_d}$ by 
\begin{equation}
\label{orthog_range_co2}
\begin{split}
    \widehat{P}_{\leq k} Z = \text{Ten}_k \left[\widehat{\bm P}_{\leq k} \bm Z^{\langle k \rangle} \right], \qquad 
     \widehat{P}_{>k} Z = \text{Ten}_k \left[ \bm Z^{\langle k \rangle} \widehat{\bm P}_{>k} \right], 
\end{split}
\end{equation}
where $\text{Ten}_k$ denotes the tensorization operator that is the inverse of the $k$th unfolding, i.e., 
$\text{Ten}_k\left(\bm Z^{\langle k \rangle}\right) = Z$. 
The projectors $\widehat{P}_{\leq j},\widehat{P}_{>k}$ commute whenever $j \leq k$ and can be used to construct the orthogonal projector onto the tangent space \cite[Corollary 3.2]{Lubich_2015} 
\begin{equation}
\label{orth_tangent_proj}
    \widehat{P}_{Y}
    = \sum_{k=1}^{d-1} 
     \widehat{P}_{\leq k-1} \widehat{P}_{>k}
    - \widehat{P}_{\leq k}  \widehat{P}_{>k}  + 
    \widehat{P}_{\leq d-1},
\end{equation}
where we set $\widehat{P}_{\leq 0}=1$.

\subsection{Oblique tangent space projectors}

The computational cost of dynamical low-rank approximation of \eqref{nonlinear_ODE0} with orthogonal tangent space projections can scale as $\mathcal{O}(n^d)$ when $G$ lacks low-rank structure. 
In this case, the orthogonal dynamical low-rank method is impractical as its computational cost is comparable to solving \eqref{nonlinear_ODE0} without low-rank compression. 
To enable efficient dynamical low-rank approximation in such cases, we introduce oblique projections onto the TT tangent space that can be computed in only $\mathcal{O}(dnr^3)$ operations for many applications where orthogonal projections require $\mathcal{O}(n^d)$. 
We construct such oblique projectors by replacing the orthogonal projectors defined in \eqref{orthog_range_co} with oblique projectors onto the same spaces 
\begin{equation}
\label{oblique_range_co}
    \begin{split}
    {\bm P}_{\leq k} = \bm U_{\leq k} \left( \bm X^{\top}_{\leq k} \bm U_{\leq k}\right)^{-1} \bm X^{\top}_{\leq k}, \qquad 
    {\bm P}_{>k} = \bm X_{>k} \left(\bm X_{>k}^{\top} \bm V_{> k}\right)^{-\top} \bm V_{>k}^{\top},
\end{split}
\end{equation}
defined for any matrices $\bm X_{\leq k} \in \mathbb{R}^{(n_1 \cdots n_k) \times r_k}$ and $\bm X_{>k} \in \mathbb{R}^{(n_{k+1} \cdots n_d) \times r_k}$ such that $\left( \bm X^{\top}_{\leq k} \bm U_{\leq k}\right)$ and $ \left(\bm X_{>k}^{\top} \bm V_{> k}\right)$ are invertible. 
Just as before, it is convenient to define oblique projectors corresponding to \eqref{oblique_range_co} that act on the tensor space $\mathbb{R}^{n_1  \times \cdots \times n_d}$ by 
\begin{equation}
\label{oblique_range_co2}
\begin{split}
    {P}_{\leq k} Z = \text{Ten}_k \left[\bm P_{\leq k} 
    \bm Z^{\langle k \rangle} \right], \qquad 
     {P}_{>k} Z = \text{Ten}_k \left[ \bm Z^{\langle k \rangle} \bm P_{>k} \right]. 
\end{split}
\end{equation}
With these projectors we can construct an oblique tangent space projector with the same form as the orthogonal tangent space projector \eqref{orth_tangent_proj}.

\begin{proposition}
\label{prop:tangent_proj}
Let $Y \in \M_{\r}$ with orthogonal decompositions of its unfolding matrices given in \eqref{svd_k} and suppose $\bm X_{\leq k}, \bm X_{>k}$ define oblique projectors \eqref{oblique_range_co}. 
Then the map 
\begin{equation}
\label{oblique_tangent_proj}
    {P}_{Y} 
    = \sum_{k=1}^{d-1} 
     {P}_{\leq k-1} {P}_{>k}
    - {P}_{\leq k} {P}_{>k}  + 
    {P}_{\leq d-1},
\end{equation}
with $P_{\leq 0} = 1$, defined for any tensor $Z \in \mathbb{R}^{n_1 \times \cdots \times n_d}$, is an oblique projector onto the tangent space $T_{Y}\M_{\r}$. 
\end{proposition}
\begin{proof}
First we show that the image of $P_{Y}$ is contained in $T_{Y} \M_{\r}$. 
Note that $T_{Y} \M_{\r}$ is a linear space and thus it is sufficient to show that the image of each term in \eqref{oblique_tangent_proj} belongs to $T_{Y} \M_{\r}$. 
To do so we utilize \cite[Corollary 3.2]{Lubich_2015} which shows that $P_{\leq j},P_{>k}$ commute whenever $j \leq k$ and 
\begin{equation}
\label{P_km1_Pk}
P_{\leq k-1} P_{>k} Z 
    = 
    \text{Ten}_{k}\left\{ \left[\bm I_{n_k} \otimes \bm P_{\leq k-1} \right] \bm Z^{\langle k \rangle} \bm P_{>k} \right\}. 
\end{equation}
Note that the result is stated for orthogonal projectors $\widehat{P}_{\leq k},\widehat{P}_{>k}$ however the proof does not require orthogonality of the projectors and thus holds for oblique projectors ${P}_{\leq k},{P}_{>k}$. 
Inserting the oblique projectors \eqref{oblique_range_co} into \eqref{P_km1_Pk} we obtain 
\begin{equation}
\label{mixed_kron}
\begin{split}
    P_{\leq k-1} P_{>k} Z 
    &= 
    \text{Ten}_{k}\left\{ \left[\bm I_{n_k} \otimes \bm U_{\leq k-1} \left( \bm X^{\top}_{\leq k-1} \bm U_{\leq k-1}\right)^{-1} \bm X^{\top}_{\leq k-1}\right] \bm Z^{\langle k \rangle} 
    \left[\bm X_{>k} \left(\bm X_{>k}^{\top} \bm V_{>k} \right)^{-T} \bm V_{>k}^{\top}\right] \right\}
    \\
    &= \text{Ten}_{k}\left\{ \left[ \bm I_{n_k} \otimes \bm U_{\leq k-1}\right] \left[\bm I_{n_k} \otimes \left( \bm X^{\top}_{\leq k-1} \bm U_{\leq k-1}\right)^{-1} \bm X^{\top}_{\leq k-1} \bm Z^{\langle k \rangle} 
    \bm X_{>k} \left(\bm X_{>k}^{\top} \bm V_{>k} \right)^{-T} \right] \bm V_{>k}^{\top} \right\} \\
    &= \text{Ten}_{k}\left\{ \left[ \bm I_{n_k} \otimes \bm U_{\leq k-1}\right] \delta \bm C_k^{\langle l \rangle} \bm V_{>k}^{\top} \right\} 
\end{split}
\end{equation}
where we used the mixed product property of the Kronecker product to obtain the second equality and defined $\delta C_k$ as the $r_{k-1} \times n_k \times r_k$ TT-core with left unfolding 
\begin{equation}
    \delta \bm C_k^{\langle l \rangle} = 
    \bm I_{n_k} \otimes \left( \bm X^{\top}_{\leq k-1} \bm U_{\leq k-1}\right)^{-1} \bm X^{\top}_{\leq k-1} \bm Z^{\langle k \rangle} \bm X_{>k} \left(\bm X_{>k}^{\top} \bm V_{>k} \right)^{-T}. 
\end{equation}
Using the identity 
\begin{equation}
\label{kron_rel}
     \bm C_{\leq k} = \left(\bm I_{n_k} \otimes \bm C_{\leq k-1} \right) \bm C_k^{\langle l \rangle},
\end{equation}
and \eqref{partial_unfolding} we can write \eqref{mixed_kron} in TT format 
\begin{equation}
\label{P+_TT}
P_{\leq k-1} P_{>k} Z  = U_{\leq k-1}
\delta C_k V_{>k}, 
\end{equation}
which has the form \eqref{tangent_element} and thus belongs to the tangent space $T_{Y} \M_{\r}$. 
For the $P_{\leq k} P_{>k}$ terms we have 
\begin{equation}
\begin{split}
    P_{\leq k}P_{>k} Z
    &= 
     \text{Ten}_{k}\left\{ \bm P_{\leq k} \bm Z^{\langle k \rangle} \bm P_{>k} \right\} \\
    &= \text{Ten}_k\left\{\bm U_{\leq k} \left[\left( \bm X^{\top}_{\leq k} \bm U_{\leq k}\right)^{-1} \bm X^{\top}_{\leq k} \bm Z^{\langle k \rangle} \bm X_{>k} \left(\bm X_{>k}^{\top} \bm V_{>k} \right)^{-T}\right] \bm V_{>k}^{\top} \right\}, 
\end{split}
\end{equation}
which we write in the TT format as 
\begin{equation}
\label{P-_TT}
P_{\leq k} P_{>k} Z = U_{\leq k} \delta \bm S_k V_{>k}, 
\end{equation}
where 
\begin{equation}
\label{ds_1}
    \delta \bm S_k =  \left( \bm X^{\top}_{\leq k} \bm U_{\leq k}\right)^{-1} \bm X^{\top}_{\leq k} \bm Z^{\langle k \rangle} \bm X_{>k} \left(\bm X_{>k}^{\top} \bm V_{>k} \right)^{-T}, 
\end{equation}
belongs to $\mathbb{R}^{r_k \times r_k}$. 
Absorbing $\delta \bm S_k$ into its left neighboring core $U_k$ we rewrite \eqref{P-_TT} as 
\begin{equation}
P_{\leq k} P_{>k} Z = U_{\leq k-1} \delta \overline{C}_k V_{>k}, 
\end{equation}
where $\delta \overline{C}_k(i_k) = U_k(i_k) \delta \bm S_k$, which has the form \eqref{tangent_element} and therefore belongs to the tangent space $T_{Y}\M_{\r}$. 

It remains to show that $P_{Y}$ is idempotent, which we verify by checking that $P_{\leq k}$ and $P_{>k}$ are idempotent for all $k = 1,\ldots,d-1$. 
For $P_{\leq k}$ we have 
\begin{equation}
\begin{split}
{P}_{\leq k}^2 Z &= P_{\leq k} \text{Ten}_k \left[\bm P_{\leq k} \bm Z^{\langle k \rangle} \right] \\
&= \text{Ten}_k\left[\bm P_{\leq k}^2 \bm Z^{\langle k \rangle} \right] \\
&= \text{Ten}_k\left[\bm P_{\leq k} \bm Z^{\langle k \rangle} \right] \\
&= {P}_{\leq k} Z, 
\end{split}
\end{equation}
for all $k=1,\ldots,d-1$. 
The idempotence of $P_{>k}$ for $k =1,\ldots,d-1$ is also easy to verify directly. 
\end{proof}

\vs
\noindent
We remark that the orthogonality of the bases $\bm U_{\leq k}$ and $\bm V_{>k}$ is not necessary for the construction of oblique tangent space projectors \eqref{oblique_tangent_proj}. 
However, imposing such orthogonality is advantageous for the numerical stability and accuracy of the projectors. 

\subsection{Interpolatory tangent space projectors}

Next we introduce a special class of oblique tangent space projectors with a cross interpolation property that enables efficient dynamical low-rank approximation. Such oblique projectors are obtained by selecting $\bm X_{\leq k}$ and $\bm X_{>k}$ in \eqref{oblique_tangent_proj} as $r_k$ columns of the identity matrix with appropriate dimension 
\begin{equation}
\label{Sdef}
\begin{split}
    \bm X_{\leq k} = \bm I_{n_1\cdots n_k}\left(:,I^{\leq k}\right), \qquad 
    \bm X_{>k} = \bm I_{n_{k+1}\cdots n_d}\left(:,I^{>k}\right).
\end{split}
\end{equation}
Here, $I^{\leq k}$ contains $r_k$ indices of the form $i_1 \cdots i_k$ and $I^{>k}$ contains $r_k$ indices of the form $i_{k+1}\cdots i_d$. 
The matrix $\bm X_{\leq k}$ can extract $r_k$ rows determined by the index sets $I^{\leq k}$ from a tensor $\bm Z^{\langle k\rangle}$ with matrix multiplication from the left 
$\bm X_{\leq k}^{\top} \bm Z^{\langle k \rangle} = \bm Z^{\langle k \rangle}(I^{\leq k}, :)$.
Similarly, the matrix $\bm X_{>k}$ can extract $r_k$ columns from $\bm Z^{\langle k\rangle}$ with matrix multiplication from the right 
$\bm Z^{\langle k \rangle} \bm X_{> k} = \bm Z^{\langle k \rangle}(:,I^{>k}).$
Moreover, when $\bm X_{\leq k}, \bm X_{>k}$ are defined as in \eqref{Sdef} the oblique projectors in \eqref{oblique_range_co} 
become interpolatory projectors 
\begin{equation}
\label{matrix_interp}
\begin{split}
\left[\bm P_{\leq k} \bm Z^{\langle k \rangle} \right]\left(I^{\leq k},:\right) = \bm Z^{\langle k \rangle}\left(I^{\leq k},:\right),
\qquad 
\left[\bm Z^{\langle k \rangle}\bm P_{>k} \right]\left(:,I^{>k}\right) = \bm Z^{\langle k \rangle}\left(:,I^{> k}\right).
\end{split}
\end{equation}
Since each index $i_1\cdots i_k$ corresponds to a multi-index $(i_1,\ldots,i_k)$ and $i_{k+1} \cdots i_d$ corresponds to a multi-index $(i_{k+1}, \ldots, i_d)$, we identify the sets of $r_k$ indices $I^{\leq k}$ and $I^{>k}$ with sets of $r_k$ multi-indices $\I^{\leq k}$ and $\I^{>k}$. 
With such identification the matrix interpolation property \eqref{matrix_interp} is equivalent to an interpolation property of the corresponding tensor operators \eqref{oblique_range_co2} 
\begin{equation}
\label{interp_prop_proj_tens}
\begin{split}
    \left[P_{\leq k} Z \right]\left( \I^{\leq k}, i_{k+1},\ldots, i_d\right) 
    = Z \left( \I^{\leq k}, i_{k+1},\ldots, i_d\right), \qquad  
    \left[ P_{>k} Z \right]\left( i_1,\ldots,i_k,\I^{>k}\right) 
    = Z \left( i_1,\ldots,i_k,\I^{>k}\right). 
\end{split}
\end{equation}
In order to obtain an oblique tangent space projectors \eqref{oblique_tangent_proj} that interpolate, we consider multi-indices that satisfy the nested conditions 
\begin{equation}
    \label{nested_inds}
    \mathcal{I}^{\leq k} \subset \mathcal{I}^{\leq k-1} \times \{1,\ldots,n_k\}, \qquad 
    \mathcal{I}^{>k} \subset \{1,\ldots,n_k\} \times \mathcal{I}^{>k+1}, \quad k = 1,2,\ldots,d-1. 
\end{equation}
To prove that \eqref{nested_inds} is sufficient for the oblique tangent space projector \eqref{oblique_tangent_proj} to interpolate, we have the following Lemma. 
\begin{lemma}
\label{lemma:P_plus_minus}
For any nested indices \eqref{nested_inds} defining oblique projectors \eqref{oblique_range_co2} and $k = 1,2,\ldots,d-1$, the projector $P_{\leq k-1} P_{>k}$ satisfies 
\begin{equation}
\left[P_{\leq k-1} P_{>k} Z \right] \left( \I^{\leq j-1},i_j,\I^{>j}\right)
= 
\begin{cases}
     Z\left(\I^{\leq j-1},i_j,\I^{>j}\right), & j = k, \\
    \left[P_{\leq k} P_{>k}  Z \right] \left( \I^{\leq j-1},i_j,\I^{>j}\right), & j > k, \\
    \left[P_{\leq k-1} P_{>k-1}  Z \right] \left( \I^{\leq j-1},i_j,\I^{>j}\right), & j < k. 
\end{cases}
\end{equation}
\end{lemma}
\begin{proof}
    The case $j = k$ follows directly from the interpolation property \eqref{interp_prop_proj_tens}. 
    For $j>k$ 
    the nesting condition \eqref{nested_inds} ensures that the first $k-1$ indices of each multi-index in $\I^{\leq j-1}$ is a multi-index in $\I^{\leq k-1}$ and that the first $k$ indices of each multi-index in $\I^{\leq j-1}$ is a multi-index in $\I^{\leq k}$. 
    Therefore we can use the interpolation property \eqref{interp_prop_proj_tens} to obtain 
    \begin{equation}
        \left[ P_{\leq k-1} Z \right]\left(\I^{\leq j-1},i_j,\ldots,i_d\right) = \left[ P_{\leq k} Z \right]\left(\I^{\leq j-1},i_j,\ldots,i_d\right) = Z\left(\I^{\leq j-1},i_j,\ldots,i_d\right). 
    \end{equation}
    Since $P_{\leq k}$ and $P_{\leq k-1}$ act only on the first $k$ dimensions of $Z$ and $P_{>k}$ acts only on dimensions $k+1,\ldots,d$ we can use the preceding equation to obtain $$\left[P_{\leq k-1} P_{>k} Z \right] \left( \I^{\leq j-1},i_j,\I^{>j}\right) = \left[P_{\leq k} P_{>k} Z \right] \left( \I^{\leq j-1},i_j,\I^{>j}\right),$$ 
    establishing the case $j>k$. 
    When $j < k$ the nested condition ensures that the $d-j$ indices of each multi-index in $\I^{>j}$ appear in multi-indices belonging to $\I^{>k}$ and $\I^{>k-1}$. 
    Therefore from the interpolation property \eqref{interp_prop_proj_tens} we have 
    \begin{equation}
        \left[ P_{>k} Z \right]\left(i_1,\ldots,i_j,\I^{>j}\right) = \left[ P_{>k-1} Z \right]\left(i_1,\ldots,i_j,\I^{>j}\right) = Z\left(i_1,\ldots,i_j,\I^{>j}\right), 
    \end{equation}
    from which we obtain the result for $j<k$. 
\end{proof}

\begin{theorem}
\label{thm:interp_proj}
For any $Y \in \M_{\r}$ and $\{\I^{\leq j},\I^{>j}\}$ nested multi-indices defining interpolatory projectors \eqref{oblique_range_co} 
the oblique tangent space projector \eqref{oblique_tangent_proj} has the cross interpolation property 
\begin{equation}
\label{tan_proj_interp}
    \left[P_{Y} Z\right] \left(\I^{\leq j-1},i_j,\I^{>j}\right) = 
     Z\left(\I^{\leq j-1},i_j,\I^{>j}\right), \qquad j = 1,2,\ldots,d. 
\end{equation}
\end{theorem}
\begin{proof}
Rearrange the terms in \eqref{oblique_tangent_proj} to write 
    \begin{equation}
        \begin{split}
            \left[P_{Y} Z\right] 
            &= \sum_{k=1}^d P_{\leq k-1} P_{>k} Z - \sum_{k=1}^{d-1} P_{\leq k} P_{>k} Z, 
        \end{split}
    \end{equation}
with $P_{>d}=1$ and evaluate at the indices $\left(\I^{\leq j-1},i_j,\I^{>j}\right)$ 
\begin{equation}
\label{interp_deriv}
\begin{split}
    \left[P_{Y} Z\right]\left(\I^{\leq j-1},i_j,\I^{>j}\right)
    &= \sum_{k=1}^d \left[P_{\leq k-1} P_{>k} Z \right] \left(\I^{\leq j-1},i_j,\I^{>j}\right)
    - \sum_{k=1}^{d-1} \left[P_{\leq k} P_{>k} Z \right]\left(\I^{\leq j-1},i_j,\I^{>j}\right). 
\end{split}
\end{equation}
Applying Lemma \ref{lemma:P_plus_minus} to each term in the first summation in \eqref{interp_deriv} we obtain 
\begin{equation}
\label{interp_deriv2}
\begin{split}
    &\sum_{k=1}^d \left[P_{\leq k-1} P_{>k} Z \right] \left(\I^{\leq j-1},i_j,\I^{>j}\right) 
    = \\
    & \qquad  
    \left(\sum_{k=1}^{j-1} \left[P_{\leq k} P_{>k} Z \right]\left(\I^{\leq j-1},:,\I^{>j}\right) \right)
    +
     Z\left(\I^{\leq j-1},i_j,\I^{>j}\right) 
    + 
    \left( \sum_{l=j+1}^{d} \left[P_{\leq l-1} P_{>l-1} Z \right]\left(\I^{\leq j-1},i_j,\I^{>j}\right) \right). 
\end{split}
\end{equation}
Re-indexing the final summation in \eqref{interp_deriv2} with $k = l-1$ and combining the result with the first summation in \eqref{interp_deriv2} yields 
\begin{equation}
\label{final_eq}
\begin{split}
    &\sum_{k=1}^d \left[P_{\leq k-1} P_{>k} Z\right] \left(\I^{\leq j-1},i_j,\I^{>j}\right)  
    =
     Z\left(\I^{\leq j-1},i_j,\I^{>j}\right) 
    + 
     \sum_{k=1}^{d-1} \left[P_{\leq k} P_{>k} Z \right]\left(\I^{\leq j-1},i_j,\I^{>j}\right). 
\end{split}
\end{equation}
Finally substituting \eqref{final_eq} into \eqref{interp_deriv} the two summations cancel and the proof is complete. 
\end{proof}

\vs
\noindent
The cross interpolation property \eqref{tan_proj_interp} of the tangent space projector resembles the interpolation property of TT-cross approximation \cite{quasi_max_vol} with nested indices. 
Next, we describe a novel index selection algorithm for constructing interpolatory projectors onto the tangent space $T_{Y}\M_{\r}$. 
%

\section{Index selection for oblique projectors and cross interpolation} 
\label{sec:index_selection}

In the preceding section we constructed tangent space projectors \eqref{oblique_tangent_proj} with the cross interpolation property \eqref{tan_proj_interp} from oblique projectors \eqref{oblique_range_co} onto the bases $\bm U_{\leq k}$ and $\bm V_{>k}$. 
We now devise an efficient algorithm based on the discrete empirical interpolation method (DEIM) for computing indices $\{I^{\leq k},I^{>k}\}$, or equivalently multi-indices $\{\I^{\leq k},\I^{>k}\}$, that yield well-defined interpolatory projectors \eqref{oblique_range_co} defined by the matrices \eqref{Sdef}. 
The DEIM (recalled in Algorithm \ref{alg:DEIM})  greedily selects indices to minimize the condition number of interpolatory projectors onto a given basis, e.g., $\bm U_{\leq k}$ or $\bm V_{>k}$, as much as possible. 
However we can not apply such algorithm directly to the matrices $\bm U_{\leq k},\bm V_{>k}$ as they have dimensions $(n_1 \cdots n_k) \times r_k$ and $(n_{k+1} \cdots n_d) \times r_k$, respectively, which is too large to store in memory and process with the DEIM algorithm. 
To address the problem of memory, we do not store the matrices $\bm U_{\leq k}$ and $\bm V_{>k}$ directly but rather the left orthogonal cores $U_j$ and right orthogonal cores $V_j$ that can be used to construct these matrices
\begin{equation}
\label{partial_prods}
\begin{split}
    \bm U_{\leq k}(i_1 \cdots i_k,:) &= U_1(i_1) \cdots U_k(i_k,:), \\
    \bm V_{>k}(i_{k+1} \cdots i_d,:) &= V_{k+1}(:,i_{k+1}) \cdots V_d(i_d), \qquad k = 1,2,\ldots,d-1. 
\end{split}
\end{equation}
To address the computational cost we propose an algorithm that only samples from a small subset of entries of the matrices $\bm U_{\leq k},\bm V_{>k}$. 
The key idea is to compute the indices $I^{\leq k}$ recursively for $k = 1,2,\ldots,d-1$ by sampling from $\bm U_{\leq k}$ only indices corresponding to multi-indices $\I^{\leq k}$ that are nested \eqref{nested_inds}. 
By considering only nested indices we reduce the number of possible indices $i_1 \cdots i_k \in I^{\leq k}$ from $n_1\cdots n_k$ to $r_{k-1} n_k$. 
We use the same idea to construct the nested index sets $\I^{>k}$ sequentially for $k=d-1,d-2,\ldots,1$. 
Since the indices obtained from this sampling approach are nested by construction, the resulting tangent space projector \eqref{oblique_tangent_proj} has the cross interpolation property \eqref{tan_proj_interp}. 
Such cross interpolation property will be useful for the dynamical low-rank approximation schemes presented in Section \ref{sec:tt_integration}. 
Hereafter we present the nested index selection algorithm in detail and then show in Theorem~\ref{thm:invertible_subs} that this nested sampling method always produces well-defined interpolatory projectors.

\subsection{The TT-cross-DEIM algorithm}
\label{sec:tt-cross-deim}

%
To compute the indices $I^{\leq j}$, begin by applying the DEIM algorithm to the $n_1 \times r_1$ matrix $\bm U_{\leq 1}$ 
\begin{equation}
\label{I_leq_1}
I^{\leq 1} = \texttt{DEIM}\left( \bm U_{\leq 1} \right). 
\end{equation} 
To obtain $I^{\leq 2}$, construct the $r_1 n_2 \times r_2$ matrix
\begin{equation}
\label{W2}
\begin{split}
\widehat{\bm U}_{\leq 2}(\alpha_{1} i_2,\alpha_2) = 
 \bm U_{\leq 2}\left(I_{\alpha_1}^{\leq 1}i_2,\alpha_2\right),
\end{split}
\end{equation}
which is the restriction of $\bm U_{\leq 2}$ to the indices $I^{\leq 1}$ with $I^{\leq 1}_{\alpha_1}$ denoting the $\alpha_1$st index in $I^{\leq 1}$. 
Then sample $r_2$ indices from this restricted matrix 
\begin{equation}
\label{l2}
    \bm l_{\leq 2} = \texttt{DEIM}\left(\widehat{\bm U}_{\leq 2}\right). 
\end{equation}
Due to the construction of $\widehat{\bm U}_{\leq 2}$ in \eqref{W2}, each index $\bm l_{\leq 2}(\alpha_2)$ corresponds to a multi-index $(\bm p_{\leq 1}(\alpha_2), \bm p_2(\alpha_2))$ where $\bm p_{\leq 1}(\alpha_2)$ identifies an index in $\I^{\leq 1}$ and $\bm p_2(\alpha_2)$ identifies an index in $\{1,\ldots,n_2\}$. 
Construct the multi-index set 
\begin{equation}
\label{I2}
    \mathcal{I}^{\leq 2}_{\alpha_2} = \left(
    \mathcal{I}^{\leq 1}_{\bm p_{\leq 1}(\alpha_2)}, \bm p_2(\alpha_2)\right), \qquad \alpha_2=1,2,\ldots,{r_2}, 
\end{equation}
which corresponds to the set of indices $I^{\leq 2}$. 
The remaining sets of indices $I^{\leq j}$ are obtained inductively with a similar procedure. 
After computing ${I}^{\leq j-1}$, construct the $r_{j-1} n_j \times r_j$ matrix 
\begin{equation}
\label{Wk}
\begin{split}
\widehat{\bm U}_{\leq j}(\alpha_{j-1} i_j,\alpha_j) = 
\bm U_{\leq j}\left(I^{\leq j-1}_{\alpha_{j-1}}i_j,\alpha_j\right), 
\end{split}
\end{equation}
which is the restriction of $\bm U_{\leq j}$ to the indices $I^{\leq j-1}$. 
Then sample $r_j$ indices from the restricted matrix 
\begin{equation}
\label{lk}
    \bm l_{\leq j} = 
    \texttt{DEIM} \left(\widehat{\bm U}_{\leq j}\right). 
\end{equation}
Due to the construction of $\widehat{\bm U}_{\leq j}$ in \eqref{Wk} each index $\bm l_{\leq j}(\alpha_j)$ corresponds to a multi-index 
$(\bm p_{\leq j-1}(\alpha_j), \bm p_j(\alpha_j))$ where $\bm p_{\leq j-1}(\alpha_j)$ identifies a multi-index in $\I^{\leq j-1}$ and $\bm p_j(\alpha_j)$ identifies an index in $\{1,\ldots,n_j\}$. 
Construct the set $\I^{\leq j}$ with multi-indices 
\begin{equation}
\label{Ik}
    \mathcal{I}^{\leq j}_{\alpha_j} = \left(
    \mathcal{I}^{\leq j-1}_{\bm p_{\leq j-1}(\alpha_j)}, \bm p_j(\alpha_j)\right), \qquad \alpha_j=1,2,\ldots,{r_j}, 
\end{equation}
which corresponds to $I^{\leq j}$. 
This procedure computes index sets $I^{\leq j}$ sequentially for $j=1,2,\ldots,d-1$ by applying the DEIM sampling algorithm to the restricted matrices $\widehat{\bm U}_{\leq j}$ of dimension $r_{j-1} n_j \times r_j$. 
In practice the restricted matrices $\widehat{\bm U}_{\leq j}$ are not obtained from $\bm U_{\leq j}$ as written in \eqref{Wk} but rather from the low-dimensional tensor cores $U_j$ that construct $\bm U_{\leq j}$ in \eqref{partial_prods}. 
%

\begin{figure}[!t]
\centering
\includegraphics[scale=0.45]{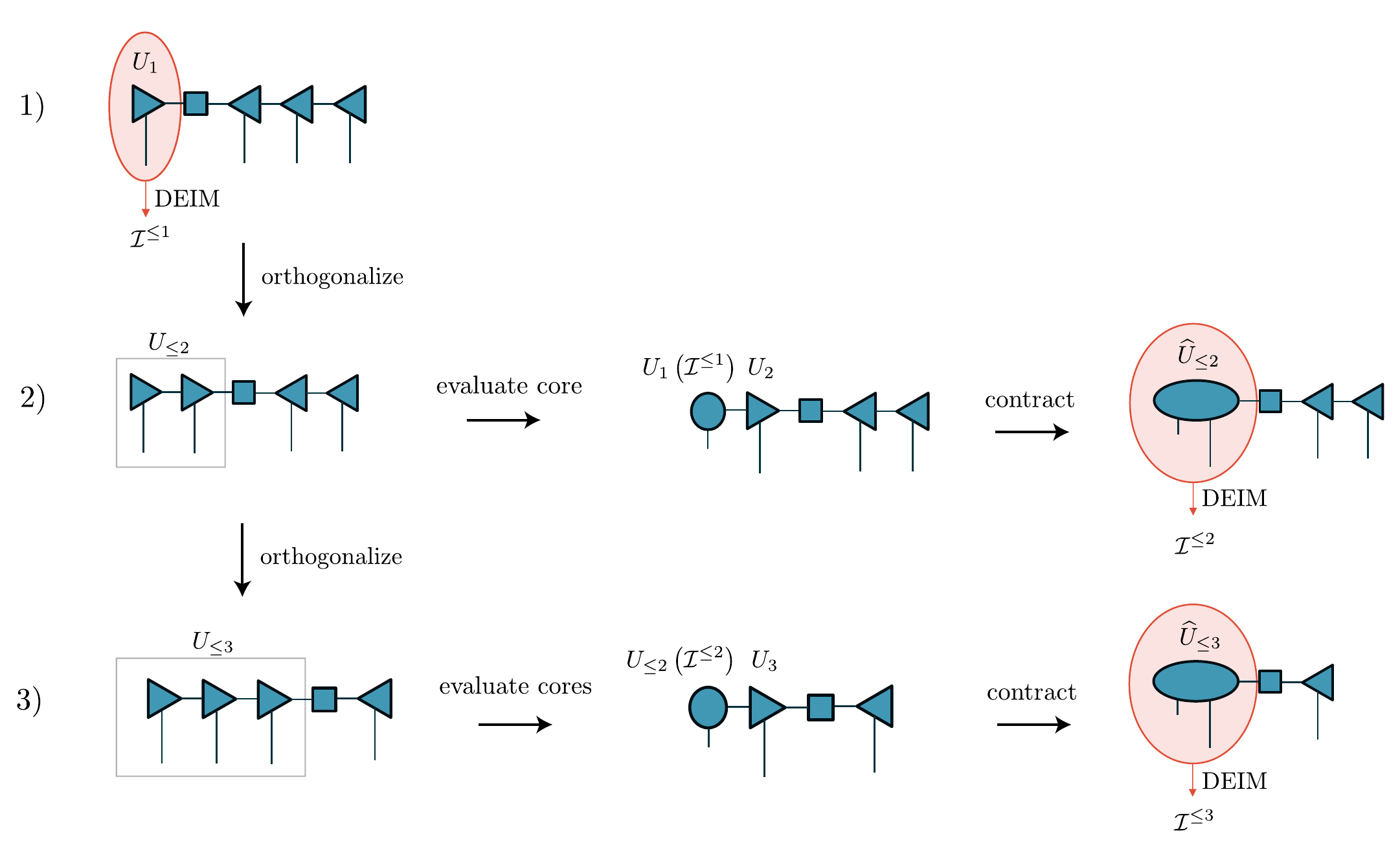}
\caption{
Illustration of the TT-cross-DEIM left-to-right sweep that computes multi-indices $\I^{\leq j}$ sequentially for $j = 1,2,\ldots,d-1$ with $d=4$. 
}
\label{fig:TT_DEIM}
\end{figure}

We compute the index sets $I^{>j}$ in a similar manner. 
First obtain $I^{>d-1}$ by sampling the $n_d \times r_{d-1}$ matrix $\bm V_{>d}$ 
\begin{equation}
\label{Idminus1}
I^{>d-1} = \texttt{DEIM} \left( \bm V_{>d-1} \right). 
\end{equation}
Then construct index sets $I^{>j}$ inductively for $j=d-2,d-3,\ldots,1$ as follows. 
After computing $I^{>j+1}$, construct the $n_{j+1} r_{j+1} \times r_j$ matrix $\widehat{\bm V}_{>j}$ 
\begin{equation}
\label{Vk}
\begin{split}
\widehat{\bm V}_{>j}\left(i_{j+1} \alpha_{j+1},\alpha_{j}\right) = 
\bm V_{>j}\left(i_{j+1} I^{>j+1}_{\alpha_{j+1}},\alpha_{j}\right) 
\end{split}
\end{equation}
which is the restriction of $\bm V_{>j}$ to the indices $\I^{>j+1}$. 
Then sample $r_j$ indices from the restricted matrix 
\begin{equation}
\label{lk_2}
\bm l_{>j} = \texttt{DEIM}\left( \widehat{\bm V}_{>j} \right). 
\end{equation}
Due to the construction of $\widehat{\bm V}_{>j}$ in \eqref{Vk} each index $\bm l_{>j}(\alpha_j)$ corresponds to a multi-index $(\bm p_{j+1}(\alpha_j),\bm p_{>{j+1}}(\alpha_{j}))$, 
where $\bm p_{j+1}(\alpha_j)$ identifies an index in $\{1,\ldots,n_{j+1}\}$ and $\bm p_{>{j+1}}(\alpha_{j})$ identifies a multi-index in $\mathcal{I}^{>j+1}$. 
Construct the set $\I^{>j}$ with multi-indices
\begin{equation}
\label{I_g_k}
    \mathcal{I}^{>j}_{\alpha_j} = \left(
    \bm p_{j+1}(\alpha_j),
    \mathcal{I}^{>j+1}_{\bm p_{>j+1}(\alpha_j)}
    \right), \qquad \alpha_j=1,2,\ldots,{r_j}, 
\end{equation}
which corresponds to $I^{>j}$. 
Just as in the computation of $I^{\leq j}$, the restricted matrices $\widehat{\bm V}_{>j}$ are not obtained from $\bm V_{> j}$ as written in \eqref{Vk} but rather from the low-dimensional tensor cores $V_j$ that construct $\bm V_{> j}$ in \eqref{partial_prods}. 

\vs
\noindent
The entire algorithm is summarized in Algorithm \ref{alg:TT_DEIM} and a tensor network diagram of the left-to-right sweep for computing $I^{\leq j}$ from the TT-cores $U_j$ is shown in Figure \ref{fig:TT_DEIM} with $d=4$. 
In Algorithm \ref{alg:TT_DEIM} we denote by $\texttt{ind2sub}$ the Matlab function that reshapes linear indices to multi-indices. 
As mentioned above, the multi-index sets ${\I^{\leq k},\I^{>k}}$ obtained in \eqref{Ik} and \eqref{I_g_k} are nested \eqref{nested_inds} by construction. 
Thus the oblique tangent space projector \eqref{oblique_tangent_proj} constructed from these indices is a cross interpolant (see Theorem \ref{thm:interp_proj}), provided it is well-defined. 
It is shown in Section \ref{sec:conditioning} that such projector is in fact well-defined. 

\begin{algorithm}
 \caption{TT-cross-DEIM index selection}
\label{alg:TT_DEIM}
\begin{algorithmic}[1]
\Require 
    \Statex $\bm U_{\leq j},\bm V_{>j}$, $j=1,\ldots,d-1$, orthogonal bases for range and co-range of $\bm Y^{\langle j \rangle}$ as in \eqref{svd_k} 
\Ensure 
    \Statex $\{\I^{\leq j},\I^{>j}\}$, 
        nested multi-index sets defining oblique projectors \eqref{oblique_range_co} via \eqref{Sdef}
\State $I^{\leq 1} = \texttt{DEIM}\left(\bm U_{\leq 1}\right)$
\Comment{left-to-right sweep: computing $I^{\leq j}$}
\For{$j = 2$ {\bf to} $d-1$ }
        \State $ \widehat{\bm U}_{\leq j} = 
        \bm U_{\leq j}\left(I^{\leq j-1} :, : \right)$ 
        \State $\bm l_{\leq j} = \texttt{DEIM}(\widehat{\bm U}_j)$
        \State $\bm p_{\leq j-1}, \bm p_j = \texttt{ind2sub}([r_{j-1},n_j],\bm l_{\leq j})$ 
        \State $\mathcal{I}^{\leq j}_{\alpha_j} = \left(\mathcal{I}^{\leq j-1}_{\bm p_{\leq j-1}(\alpha_j)}, \bm p_j(\alpha_j) \right), \qquad \alpha_j=1,2,\ldots,r_j$
\EndFor

\State $I^{>d-1} 
= \texttt{DEIM}\left(\bm V_{>d-1}\right)$
\Comment{right-to-left sweep: computing $I^{>j}$}
\For{$j = d-2$ {\bf to} $1$ }
    \State $\widehat{\bm V}_{>j} = 
    \bm V_{>j}\left(: I^{>j+1}, : \right)$
   \State $\bm l_{>j} = \texttt{DEIM}\left( \widehat{\bm V}_{>j} \right)$ 
    \State $\bm p_{j+1}, \bm p_{>j+1} = \texttt{ind2sub}([n_j, r_{j}],\bm l_{>j})$ 
    \State $\mathcal{I}^{>j}_{\alpha_{j}} = 
    \left(\bm p_{j+1}(\alpha_{j}), 
    \mathcal{I}^{>j+1}_{\bm p_{>j+1}(\alpha_{j})}\right), \qquad \alpha_{j} =1,2,\ldots,r_{j}$ 
\EndFor
\end{algorithmic}
\end{algorithm}

\paragraph{Computational cost}
%
To simplify the operation count of the proposed TT-cross-DEIM algorithm, we assume that $r_k = r$ and $n_k = n$ for all $k=1,2,\ldots,d$. 
The Algorithm requires access to all $d-1$ orthogonal representations \eqref{orth_TT} which can be computed in $\mathcal{O}(dnr^3)$ operations \cite{OseledetsTT}. 
With these orthogonal representations available, each index set $I^{\leq k}$ and $I^{>k}$ for $k=1,2,\ldots,d-1$ is computed by applying DEIM to a matrix of size $rn \times r$ each of which requires $\mathcal{O}(nr)$ operations. 
Therefore the total number of operations in the TT-cross-DEIM algorithm is $\mathcal{O}(dnr^3)$. 

\vs
\noindent
Note that the operation count of TT-cross-DEIM is dominated by the computation of orthogonal TT representations and thus has the same complexity as many common TT algorithms, e.g., TT rounding. 

\subsection{Condition of the oblique projectors}
\label{sec:conditioning}

For the oblique projectors \eqref{oblique_range_co} to be defined (and thus for the oblique tangent space projector \eqref{oblique_tangent_proj} to be defined) the $r_j \times r_j$ matrices $\bm M_j = \bm X_{\leq j}^{\top} \bm U_{\leq j}$ and $\bm N_j = \bm X_{>j}^{\top} \bm V_{>j}$
must be invertible. 
For interpolatory projectors, i.e., when $\bm X_{\leq j}$ and $\bm X_{>j}$ are defined as in \eqref{Sdef}, these matrices are given entry-wise by 
\begin{equation}
\label{Mk}
\bm M_j(\alpha_j,\beta_j) = 
\bm U_{\leq j}\left(I^{\leq j}_{\alpha_j},\beta_j\right),
\qquad 
\bm N_j(\alpha_j,\beta_j) = 
\bm V_{>j}\left(I^{>j}_{\alpha_j},\beta_j \right). 
\end{equation}
The following result shows that the TT-cross-DEIM produces indices that yield invertible matrices \eqref{Mk} and therefore define oblique projectors \eqref{oblique_range_co} and \eqref{oblique_tangent_proj}. 
\begin{theorem}
\label{thm:invertible_subs}
If $Y \in \M_{\r}$ and $I^{\leq j},I^{>j}$ are obtained with the TT-cross-DEIM 
then the $r_j \times r_j$ matrices 
\eqref{Mk} are invertible for all $j = 1,2,\ldots,d-1$. 
\end{theorem}
\begin{proof}
We prove the result for $\bm M_j$ by induction on $j$. 
The $r_1$ indices $I^{\leq 1}$ are obtained in \eqref{I_leq_1} from $\bm U_{\leq 1}$ with the DEIM. 
Because $Y \in \M_{\r}$, the unfolding matrix $\bm U_{\leq 1}$ is full rank. 
Thus we have from \cite[Lemma 3.1]{DEIM-CUR} that the $r_1 \times r_1$ matrix $\bm M_1 = \bm U_{\leq 1}\left(I^{\leq 1},:\right)$ is full rank, establishing the result for $\bm M_j$ when $j=1$. 
Now assume that $\bm M_{j-1}$ is full rank. 
Rewriting \eqref{Wk} using \eqref{partial_tensor}-\eqref{partial_unfolding} we have 
\begin{equation}
\begin{split}
\widehat{\bm U}_j(\alpha_{j-1} i_j,\alpha_j) &= 
\bm U_{\leq j-1}\left(I^{\leq j-1}_{\alpha_{j-1}},:\right) 
\bm U_j^{\langle r \rangle}(:,i_j \alpha_j) 
\\
&= \bm M_{j-1}(\alpha_{j-1},:) \bm U_j^{\langle r \rangle}(:,i_j \alpha_j). 
\end{split}
\end{equation}
By assumption $\bm M_{j-1}$ is full rank and since $Y \in \M_{\r}$ both unfoldings of $\bm U_j$ are full rank. 
It follows that $\widehat{\bm U}_j$ is full rank. 
In \eqref{lk} the $r_j$ indices $\bm l_{\leq j}$ are obtained from $\widehat{\bm U}_j$ with the DEIM and invoking \cite[Lemma 3.1]{DEIM-CUR} we have that the $r_j \times r_j$ matrix $\widehat{\bm U}_j(\bm l_{\leq j},:)$ is full rank. 
The indices $\I^{\leq j}$ are obtained from $\bm l_{\leq j}$ in \eqref{Ik} so that 
\begin{equation}
\bm U_{\leq j}\left(I^{\leq j}_{\alpha_j},\beta_j\right) = \widehat{\bm U}_{\leq j}\left(\bm l_{\leq j}(\alpha_j),\beta_j\right),
\end{equation}
proving the result for $\bm M_j$.  
The statement for $\bm N_j$ is proven similarly with induction on $j = d-1,d-2,\ldots,1$. 
\end{proof}

\vs 
\noindent

\vs
\noindent
In addition to generating invertible matrices \eqref{Mk}, the TT-cross-DEIM is a greedy algorithm that aims to minimize the condition number of these matrices as much as possible while ensuring the indices remain nested. 
Indeed, the DEIM algorithm is used to compute $\bm l_{\leq j}$ in \eqref{lk} and $\bm l_{>j}$ in \eqref{lk_2}, selecting indices greedily to keep the condition numbers of $\widehat{\bm U}_{\leq j}(\bm l{\leq j},:)$ and $\widehat{\bm V}_{>j}(\bm l{>j},:)$ small. 
Moreover, $\widehat{\bm U}_{\leq j}$ is the restriction of $\bm U_{\leq j}$ to the indices $I^{\leq j-1}$ (see \eqref{Wk}), and $\widehat{\bm V}_{>j}$ is the restriction of $\bm V_{>j}$ to the indices $I^{>j+1}$ (see \eqref{Vk}). 
We also note that other sparse sampling methods can be used in place of DEIM in the TT-cross-DEIM algorithm, e.g., Q-DEIM or oversampling methods, which can yield better conditioned interpolatory projectors in some cases.

\subsection{Tensor cross interpolation}

We have shown above that the multi-indices obtained with the TT-cross-DEIM produce a well-defined interpolatory projector \eqref{oblique_tangent_proj} onto the tangent space. 
Incidentally, we can use the same multi-index sets to parameterize $Y$ with tensor cross interpolation. 
Recall that TT-cross approximation \cite{parr_TT_cross_Dolgov,TT-cross,quasi_max_vol,
cross_error} is a specific instance of the TT format \eqref{TT_SVD} that generalizes the matrix CUR decomposition to tensors. 
In this representation the TT-cores are defined by the entries of $Y$ 
\begin{equation}
\label{TT_cross}
    \tilde{Y}(i_1,\ldots,i_d) 
    = 
    \prod_{k=1}^{d-1} Y\left(\I^{\leq k-1},i_k,\I^{>k}\right) \left[Y\left(\I^{\leq k},\I^{>k}\right)\right]^{-1} 
    Y\left(\I^{\leq d-1},i_d\right), 
\end{equation}
where for convenience we set $\I^{\leq 0} = \emptyset$. 
It is well-known that the nested condition \eqref{nested_inds} is sufficient for the tensor cross approximation \eqref{TT_cross} to be a tensor cross interpolant \cite{quasi_max_vol} 
\begin{equation}
\label{interp_property}
\tilde{Y}\left(\mathcal{I}^{\leq k-1}, i_k, \mathcal{I}^{>k}\right) =  
Y\left(\mathcal{I}^{\leq k-1},i_k, \mathcal{I}^{>k}\right), \qquad k = 1,2,\ldots,d. 
\end{equation}
We now use the nested multi-index sets constructed by the TT-cross-DEIM to prove that any TT can be exactly represented as a TT-cross interpolant. 
This result follows as a Corollary of Theorem \ref{thm:invertible_subs}. 

\begin{corollary}
\label{corr:exact_cross}
Any $Y \in \M_{\r}$ can be exactly represented as a rank-$\bm r$ TT-cross interpolant with nested indices. 
\end{corollary}
\begin{proof}
Let $\{I^{\leq j},I^{>j}\}$ be nested index sets obtained with the TT-cross-DEIM. 
Using \eqref{svd_k} write the $r_j \times r_j$ matrix $Y\left(\I^{\leq j},\I^{>j}\right)$ as 
\begin{equation}
\label{MN}
\begin{split}
Y\left(\I^{\leq j},\I^{>j}\right) &= \bm U_{\leq j}\left(I^{\leq j},:\right) \bm S_j \bm V_{>j}\left(I^{>j},:\right)^{\top} 
\\ 
&= \bm M_j \bm S_j \bm N_j^{\top}, 
\end{split}
\end{equation}
where we used \eqref{Mk} to obtain the second equality. 
We have shown in Theorem \ref{thm:invertible_subs} that $\bm M_j$ and $\bm N_j$ are invertible and since $Y \in \M_{\r}$, the matrices $\bm S_j$ are also invertible for all $j=1,2,\ldots,d-1$. 
Therefore the matrices $Y\left(\I^{\leq j},\I^{>j}\right)$ are invertible for all $j=1,2,\ldots,d-1$ and hence by \cite[Theorem 2]{parr_TT_cross_Dolgov} the nested multi-indices $\{\I^{\leq j},\I^{>k}\}$ provide an exact representation of $Y$ with TT-cross interpolation. 
\end{proof}

\vs
\noindent
Several algorithms for computing tensor cross approximations \eqref{TT_cross} from black-box tensors based on the maximum volume principle have recently been developed \cite{TT-cross,parr_TT_cross_Dolgov}. 
The purpose of our proposed TT-cross-DEIM algorithm is to obtain interpolatory tangent space projections for a tensor $Y \in \M_{\r}$, not for black-box tensor approximation. 
However, it was recently demonstrated in \cite{ghahremani2024cross} that DEIM-based cross approximation algorithms can be applied iteratively to obtain tensor cross approximations from black-box tensors with comparable performance to the corresponding maximum volume algorithms. 

\section{Time integration on tensor train manifolds}
\label{sec:tt_integration}

We now consider the dynamical low-rank evolution equation \eqref{DLRA} for tensors ($d \geq 2$) using interpolatory projectors \eqref{oblique_tangent_proj} onto tangent spaces of TT manifolds. 
The concept of dynamical low-rank tensor approximation is a natural extension the dynamical low-rank matrix approximation described in Section~\ref{sec:matrix_DLR}. Similar to the matrix case, classical dynamical low-rank tensor approximation uses the orthogonal projector \eqref{orth_tangent_proj} to obtain the best approximation (in the Frobenius norm) of $G(Y,t)$ in the tangent space of the TT manifold (see Figure \ref{fig:tensor_mfld}). However, as noted earlier, orthogonal projection onto the TT tangent space can have a computational cost $\mathcal{O}(n^d)$ when $G$ lacks low-rank structure.
By replacing the orthogonal projector with an interpolatory projector onto the TT tangent space we propose new dynamical low-rank methods with computational cost scaling as $\mathcal{O}(dnr^3)$ for a large class of nonlinear functions $G$ that do not have rank structure. 
In particular, the proposed interpolatory dynamical low-rank tensor methods are efficient whenever it is possible to evaluate the tensor $G(X,t)$ entry-wise. 

We also point out that cross approximation algorithms based on the maximum volume principle developed in \cite{TT-cross} are designed to obtain TT approximations of tensors that can be evaluated entry-wise. TT-cross based on maximum volume can be used to obtain a low-rank approximation of the tensor $G(Y,t)$ at each time step. Such approximation can then be projected orthogonally onto the tangent space for a dynamical low-rank method or used in a step-truncation scheme. The TT-cross-DEIM index selection strategy developed in the present work differs in that it selects interpolation indices from the solution tensor $Y(t)$, not from $G(Y,t)$. Such indices are selected so that $Y$ can be represented using a TT-cross interpolant. More importantly, it allows $G(Y,t)$ to be interpolated directly onto the tangent space of the TT manifold at $Y$, making the TT-cross-DEIM particularly suitable for efficient dynamical low-rank approximation.

%
Hereafter we propose two time integration schemes for solving the dynamical low-rank equation \eqref{DLRA} with interpolatory TT tangent space projectors \eqref{oblique_tangent_proj}. 
The first scheme, referred to as TT-cross time integration, extends the matrix cross integrator described in Section~\ref{sec:matrix_cross} to tensors in the TT format. This method integrates forward in time the entries of the solution tensor $Y(t)$ required to construct the TT-cross interpolant \eqref{TT_cross} at any time $t$. 
The second scheme extends the projector-splitting scheme described in Section~\ref{sec:proj_split_mat} to tensors in the TT format. It is a direct generalization of the projector-splitting integrator for orthogonal dynamical low-rank approximation introduced in \cite{lubich2014projector} for matrices and subsequently generalized to TTs \cite{Lubich_2015}, Tucker tensors \cite{LubichTucker} and tree tensor networks \cite{Lubich_tree_tensor}, to interpolatory tangent space projectors in the TT format. 

\subsection{TT-cross integrator}
\label{sec:TT_cross_int}
The time-dependent interpolatory TT tangent space projector \eqref{oblique_tangent_proj} in \eqref{DLRA} is defined at each time $t$ by a set of time-dependent multi-indices $\{\I^{\leq k}(t),\I^{>k}(t)\}$. 
Selecting such with the TT-cross-DEIM ensures that they are nested \eqref{nested_inds} at each time $t$. Hence the tangent space projector has the cross interpolation property \eqref{tan_proj_interp} at each $t$. 
Evaluating \eqref{DLRA} at the multi-indices $\{\I^{\leq k}(t),\I^{>k}(t)\}$ and utilizing the cross interpolation property yields evolution equations for the entries of $Y(t)$ defining a TT-cross interpolant 
\begin{equation}
\label{propagator_interp}
\displaystyle\frac{d Y\left(\I^{\leq k-1}(t),:,\I^{>k}(t),t\right)}{dt} = 
G_Y\left(\I^{\leq k-1}(t),:,\I^{>k}(t),t\right), \vspace{.2cm} \qquad k =1,2,\ldots,d, 
\end{equation} 
where we defined the tensor 
$G_Y(t) = G\left(Y(t), t \right)$.  
Equation \eqref{propagator_interp} consists of $\sum_{k=1}^d r_{k-1} n_k r_k$ coupled nonlinear differential equations governing the evolution of a subset of entries in the approximate solution $Y(t) \in \M_{\r}$, which can be integrated using standard methods. 
If $G$ arises from the spatial discretization of a PDE \eqref{nonlinear_PDE0} involving differential operators then evaluating $G_Y\left(\I^{\leq k-1},:,\I^{>k}\right)$ requires entries of $Y$ at indices adjacent to $\{\I^{\leq k},\I^{>k}\}$. 
Letting $\I^{\leq k}_{(a)}$, $\I^{>k}_{(a)}$ denote the union of $\I^{\leq k}$, $\I^{>k}$ with the required adjacent indices, the right-hand side tensors in \eqref{propagator_interp} are computed by 
\begin{equation}
\label{RHS_G}
G_Y\left(\I^{\leq k-1},:,\I^{>k}, t \right) = 
G\left(Y\left(\I^{\leq k-1}_{(a)},:,\I^{>k}_{(a)}, t \right), t \right), \qquad k = 1,2,\ldots,d. 
\end{equation}
The values of $Y$ at adjacent indices can always be obtained by constructing the low-rank solution $Y(t)$ in the TT format using TT cross interpolation as described below. 
Computing \eqref{RHS_G} is efficient for any nonlinear $G$ that we can evaluate entry-wise, regardless of the low-rank structure in $G$. 
This gives the evolution equations \eqref{propagator_interp} a clear computational advantage over the evolution equations of orthogonal dynamical low-rank approximation or step-truncation methods \cite{rodgers2020step-truncation}, which require a low-rank representation of $G$ to be practical. 
We also note that the stiffness of the evolution equations \eqref{propagator_interp} is independent of the singular values of the solution tensor $Y(t)$, unlike other dynamical low-rank methods. However, the stability of constructing the solution in TT format using cross interpolation, which is often needed to evaluate the right-hand side of \eqref{propagator_interp}, depends on the condition of the interpolatory projectors obtained through the TT-cross-DEIM index selection as shown below. 

\subsubsection{Constructing the low-rank solution in TT format}
\label{sec:constructing_solution}

We can access entries of the approximate solution $Y(t)$ at indices other than the interpolation indices by constructing $Y(t)$ using TT-cross interpolation \eqref{TT_cross}. 
\begin{equation}
\label{TT_cross_time}
Y(i_1,\ldots,i_d,t) = 
\prod_{k=1}^{d-1} 
Y\left(\I^{\leq k-1}(t),i_k, \I^{>k}(t),t\right) 
\left[ Y\left(\I^{\leq k}(t), \I^{>k}(t),t\right) \right]^{-1} 
Y\left(\I^{\leq d-1}(t),i_d,t\right), 
\end{equation}
Such entries are often needed to evaluate the right-hand side of \eqref{propagator_interp} and the TT-representation \eqref{TT_cross_time} is also needed to construct indices for interpolatory projection onto the tangent space. 
Constructing $Y(t)$ using \eqref{TT_cross_time} can lead to numerical instability as the time-dependent $r_k \times r_k$ matrices $Y\left(\I^{\leq k},\I^{>k},t\right)$ can be ill-conditioned. 
Hereafter we describe a more robust method for computing $Y(t)$ by orthogonalization, omitting the dependence on $t$ to simplify notation. 
Take QR-decompositions 
\begin{equation}
\label{fk_orth}
\left[\bm Y\left(\I^{\leq k-1},:,\I^{>k}\right)\right]^{\langle l \rangle} = 
\bm Q_k^{\langle l \rangle} \bm R_k, \qquad k=1,2,\ldots,d-1,
\end{equation}
to write 
\begin{equation}
\label{fk_orth2}
Y\left(\I^{\leq k-1},i_k,\I^{>k}\right) = 
Q_k(i_k) \bm R_k, \qquad 
Y\left(\I^{\leq k},\I^{>k}\right) = 
\bm Q_k^{\langle l \rangle}\left(\bm l_{\leq k}, : \right) \bm R_k,
\end{equation}
where $\bm l_{\leq k}$ is defined in \eqref{lk}. 
Then substitute \eqref{fk_orth2} into \eqref{TT_cross_time} to obtain 
\begin{equation}
\label{TT_cross_time2}
\begin{split}
Y(i_1,\ldots,i_d) 
&= \prod_{k=1}^{d-1} 
Q_k(i_k) \bm R_k \left[\bm Q_k^{\langle l\rangle} \left(\bm l_{\leq k}, : \right) \bm R_k \right]^{-1} Y\left(\I^{\leq d-1},i_d\right) \\ 
&= \prod_{k=1}^{d-1} 
Q_k(i_k) \left[\bm Q_k^{\langle l\rangle} \left(\bm l_{\leq k}, : \right) \right]^{-1} Y\left(\I^{\leq d-1},i_d\right),
\end{split}
\end{equation}
Computing $Y$ via \eqref{TT_cross_time2} instead of \eqref{TT_cross_time} yields a more stable numerical algorithm as the matrices $\widehat{\bm Q}_k = \bm Q_k^{\langle l\rangle} \left(\bm l_{\leq k}, : \right)$ are related to the orthogonal bases $\bm U_{\leq k}$ from which the multi-index sets $\I^{\leq k},\I^{>k}$ were obtained with the TT-cross-DEIM index selection algorithm and therefore are have smaller condition number than $Y\left(\I^{\leq k},\I^{>k}\right)$. 
The improvement in condition number is verified by our numerical experiments as shown in Figure \ref{fig:AC3D}(d).

\subsubsection{Discrete-time TT-cross integrator}
Let us describe one step of the TT-cross time integration scheme from time $t_0$ to $t_1 = t_0+\Delta t$ starting from the rank-$\bm r$ TT representation 
\begin{equation}
\label{TT_time}
Y(t_0) = C_1(t_0) C_2(t_0) \cdots C_d(t_0). 
\end{equation}
First compute the indices $\left\{\I^{\leq k}(t_0),\I^{>k}(t_0) \right\}$ for the interpolatory projector defining the dynamical low-rank evolution equation \eqref{DLRA} using the TT-cross-DEIM algorithm. 
Then integrate the evolution equations \eqref{propagator_interp} from time $t_0$ to $t_1$ using an explicit time stepping scheme with multi-indices fixed at time $t_0$, e.g., Euler forward yields 
\begin{equation}
\label{advanced_slices}
\begin{split}
Y\left(\I^{\leq k-1}(t_0),:,\I^{>k}(t_0) , t_1 \right) = 
Y\left(\I^{\leq k-1}(t_0),:,\I^{>k}(t_0) , t_0\right) 
 + \Delta t 
 G_Y\left(\I^{\leq k-1}(t_0),:,\I^{>k}(t_0) , t_0\right),
\end{split}
\end{equation}
for all $k = 1,2,\ldots,d$. 
Use the result of explicit time integration \eqref{advanced_slices} to construct TT-cores for the solution at time $t_1$ 
\begin{equation}
\label{cores_t_next}
Y(t_1) = C_1(t_1) C_2(t_1) \cdots C_d(t_1), 
\end{equation}
with the QR-stabilized procedure described in \eqref{fk_orth}-\eqref{TT_cross_time2}, i.e., 
\begin{equation}
\label{updated_cores}
\begin{split}
&C_k(i_k,t_1) = Q_k(i_k,t_1) \left[\bm Q_k^{\langle l\rangle} \left(\bm l_{\leq k}(t_0), :,t_1\right) \right]^{-1}, \quad k =1, 2,\ldots,d-1, \\
&C_d(t_1) = Y\left(\I^{\leq d-1}(t_0),:,t_1\right). 
\end{split}
\end{equation}
This completes one step of the TT-cross time integration scheme.

\paragraph{Computational cost}

To simplify the operation count of one step of the TT-cross integrator, we assume that $r_k = r$ and $n_k = n$ for all $k=1,2,\ldots,d$. 
As shown in Section \ref{sec:index_selection}, the TT-cross-DEIM algorithm used to obtain the indices $\{\I^{\leq k},\I^{>k}\}$ requires $\mathcal{O}(dnr^3)$ operations. 
Preparing the tensors $Y(\I^{\leq k-1},:\I^{>k})$ required to take the explicit time step \eqref{advanced_slices} requires $d-1$ matrix multiplications with matrices of size $r \times r$ and $r \times (nr)$ for a number of operations scaling as $\mathcal{O}(dnr^3)$. 
For many $G$ that can be evaluated entry-wise (e.g., entry-wise nonlinearities), the cost of evaluating $G_Y(^{\leq k-1},:,\I^{>k};t)$ is on the same order of computing the subtensors $Y(\I^{\leq k-1},:\I^{>k})$, i.e., $\mathcal{O}(dnr^3)$. 
Finally reconstructing the tensor cores at time $t_1$ in \eqref{updated_cores} requires $d-1$ QR-decompositions of matrices with size $r \times (nr)$ and inverting $d-1$ matrices of size $r \times r$ with total cost scaling as $\mathcal{O}(dnr^3)$. 
Thus the computational cost of one step of the TT-cross time integrator scales as $\mathcal{O}(dnr^3)$. 
We note that it is not strictly necessary to perform these steps at every time step as multi-indices can be reused over many time steps, provided the condition of $\bm Q_k^{\langle l\rangle} \left(\bm l_{\leq k}, :,t\right)$ remains under control. 
If the computation of new indices is not required then \eqref{advanced_slices} can be iterated for a number of time steps at a cost of $\mathcal{O}(dnr^2)$ operations before computing new indices.

\subsection{Interpolatory projector-splitting integrator}
\label{sec:proj_split}

Next, we propose a second time integration scheme by directly applying a splitting integrator to the dynamical low-rank evolution equation \eqref{DLRA}. This method is a direct generalization of the orthogonal projector-splitting integrator introduced in \cite{Lubich_2015} for TTs. 
As we will see, the oblique projector-splitting integrator satisfies the same desirable properties: it is robust to small singular values, it exactly reproduces low-rank solutions, and one step of the integrator can be implemented as an efficient sweeping algorithm. 
The derivation of the integrator and the proofs of these results follow the same steps as the orthogonal projector-splitting integrator. 
Inserting the oblique tangent space projector \eqref{oblique_tangent_proj} into the dynamical low-rank evolution equation \eqref{DLRA} we see that the right-hand side is a sum of $2d-1$ terms 
\begin{equation}
    P_Y G(Y,t) = \sum_{j=1}^{d-1} P_j^+ G(Y,t) - P_j^-G(Y,t) + P_d^+ G(Y,t), 
\end{equation}
where $P_j^+ = P_{\leq j-1}P_{>j}$ and $P_j^- = P_{\leq j}P_{>j}$. 
Integrating \eqref{DLRA} from time $t_0$ to $t_1 = t_0 + \Delta t$ with first order Lie-Trotter splitting requires solving the $2d-1$ substeps 
\begin{equation}
\label{split_diffs}
    \begin{aligned}
        \frac{d Y_1^+(t)}{dt} 
        &= P_1^+ G\left(Y_1^+,t\right), 
        && Y_1^+(t_0) = Y(t_0), 
        \\
        \frac{d Y_1^-(t)}{dt} 
        &= -P_1^- G\left(Y_1^-,t\right), 
        \quad && Y_1^-(t_0) = Y_1^+(t_1), 
        \\
        & \ \ 
        \vdots 
        \\
        \frac{d Y_j^+(t)}{dt} 
        &= P_j^+ G\left(Y_j^+,t\right), 
        && Y_j^+(t_0) = Y_{j-1}^-(t_1), 
        \\
        \frac{d Y_j^-(t)}{dt} 
        &=  -P_j^- G\left(Y_j^-,t\right), 
        && Y_j^-(t_0) = Y_j^+(t_1),
        \\
         & \ \ 
         \vdots 
         \\
         \frac{d Y_{d}^+(t)}{dt} 
         &= P_d^+ G\left(Y_d^+,t\right), 
         && Y_d^+(t_0) = Y_{d-1}^-(t_1), 
    \end{aligned}
\end{equation}
in consecutive order to obtain the approximate solution $Y(t_1) = Y_d^+(t_1)$ at time $t_1$. 
Note that the projectors $P_j^+, P_j^-$ depend on solutions to each substep $Y_j^+(t)$ or $Y_j^-(t)$ and thus are time-dependent. 
In the case of orthogonal tangent space projector-splitting it was shown in \cite[Theorem 4.1]{Lubich_2015} that $P_j^+, P_j^-$ can be kept constant during each substep and each of the differential equations \eqref{split_diffs} can be solved exactly by updating a single TT-core. 
We have an analogous result for the oblique projector-splitting integrator. 
In the following Theorem we suppress the dependence of the multi-indices $\I^{\leq k},\I^{>k}$ on $t$ although it is assumed that such indices are selected at each time $t$ so that the interpolatory tangent space projector \eqref{oblique_tangent_proj} is well-defined. 

\begin{theorem}
\label{thm:one_core_updates}

Each split differential equation in \eqref{split_diffs} is solved exactly using time-independent projectors $P_j^+$ and $P_j^-$ at $Y_j^+(t_0)$ and $Y_j^-(t_0)$, respectively. 
%
%
Moreover, if $Y_j^+(t_0)$ has the TT representation 
    \begin{equation}
    \label{Yj+0}
     Y_j^+(t_0) = U_{\leq j-1} \left[ U_j \bm S_j \right] V_{>j}
     \end{equation}
    then 
    $$Y_j^+(t) = U_{\leq j-1} K_j(t) V_{>j},$$ 
    where 
    \begin{equation}
    \label{one_core_plus}
        \frac{d K_j(t)}{dt} = 
        \left[U_{\leq j-1}\left(\I^{\leq j-1},:\right)\right]^{-1} 
        G_j^+\left(\I^{\leq j-1},:,\I^{>j},t\right) 
    \left[V_{>j}\left(:,\I^{>j}\right)\right]^{-\top}, \quad K_j(t_0) = U_j \bm S_j,
    \end{equation}
    and $G_j^+(t) = G\left(Y_j^+(t),t\right)$. 

    \noindent
    Similarly if $Y_j^-(t_0)$ has the TT representation 
    $
    Y_j^-(t_0) = U_{\leq j} \bm S_j(t_0) V_{>j} 
    $
    then 
    $$
    Y_j^-(t) = U_{\leq j} \bm S_j(t) V_{>j}, 
    $$ 
    where 
    \begin{equation}
    \label{one_core_minus}
        \frac{d \bm S_j(t)}{dt} 
        = 
        -\left[U_{\leq j}\left(\I^{\leq j},:\right)\right]^{-1} 
        G_j^-\left(\I^{\leq j},\I^{>j},t\right) 
        \left[V_{>j}\left(:,\I^{>j}\right)\right]^{-\top},
    \end{equation}
    and $G_j^- = G\left(Y_j^-(t),t\right)$. 
\end{theorem}
\begin{proof}
The proof follows a similar approach to the analogous proof for the orthogonal projector-splitting integrator. 
First recall that we have shown in the proof of Proposition \ref{prop:tangent_proj} that $P_{\leq j-1}P_{>j}$ maps onto a tangent space of $\M_{\r}$ at each time $t$. 
This ensures that $Y_j^+(t)$ belongs to $\M_{\r}$ for all $t$ and therefore admits a time-dependent orthogonalized rank-$\bm r$ TT decomposition of the form 
\begin{equation}
\label{fkplus_svd}
    Y_j^+(t) = U_{\leq j-1}(t) K_j(t) V_{>j}(t). 
\end{equation}
Substituting \eqref{fkplus_svd} into \eqref{split_diffs} and using the product rule we obtain 
\begin{equation}
\label{prod_+}
\begin{split}
&\frac{d U_{\leq j-1}(t)}{dt} K_j(t) V_{>j}(t) 
 + 
 U_{\leq j-1}(t) \frac{ d K_j(t)}{dt} V_{>j}(t) 
 + U_{\leq j-1}(t) K_j(t) \frac{d V_{>j}(t)}{dt} \\
 & \quad = P_{\leq j-1} P_{>j} G\left(Y_j^+(t),t\right) \\
 & \quad = U_{\leq j-1}(t) \delta C_j(t) V_{>j}(t), 
\end{split}
\end{equation}
where we used \eqref{P+_TT} to obtain the third line with  
\begin{equation}
\begin{split}
    \delta C_j(t) = 
    \left[U_{\leq j-1}\left(\I^{\leq j-1},:,t\right)\right]^{-1} 
     G_j^+\left(\I^{\leq j-1},:,\I^{>j},t\right) 
    \left[V_{>j}\left(:,\I^{>j},t\right)\right]^{-\top}
\end{split}
\end{equation}
Equation \eqref{prod_+} is solved exactly by setting 
$d U_{\leq j-1}(t)/dt =0$, 
$d V_{>j}(t)/dt=0$ and 
$d K_j(t)/dt = \delta C_j(t)$ and from the initial condition $Y_j^+(t_0)$ in \eqref{Yj+0} we obtain 
\begin{equation}
\begin{split}
&U_{\leq j-1}(t) = U_{\leq j-1}, \quad V_{>j}(t) = V_{>j}, \quad K_j(t_0) = U_j \bm S_j, 
\end{split}
\end{equation}
proving the result for $Y_j^+(t)$. 
%
%
%
The proof of the assertion for $Y_j^-(t)$ is similar. 
\end{proof}

\vs
\noindent
Similar to the TT-cross evolution equations \eqref{propagator_interp}, computing the right-hand side of the differential equation \eqref{one_core_plus} requires evaluating $G_j^+$ at a subset of $r_{j-1} n_j r_j$ indices and computing the right-hand side of \eqref{one_core_minus} requires evaluating $G_j^-$ at a subset of $r_j^2$ indices. 
These evaluations are efficient for any $G$ that can be evaluated entry-wise and do not require $G$ to have any low-rank structure. 
The differential equations \eqref{one_core_plus} and \eqref{one_core_minus} involve inverses of $r_j \times r_j$ matrices $U_{\leq j}\left(\I^{\leq j},:\right)$ and $V_{>j}\left(:,\I^{>j}\right)$. 
These matrices define the interpolatory projectors \eqref{oblique_range_co} and we select the multi-indices $\I^{\leq j},\I^{>j}$ with the TT-cross-DEIM at each time $t$ to keep their condition number is small during time integration. 

\subsubsection{Sweeping algorithm for interpolatory projector-splitting integrator}
\label{sec:sweep}

One complete step of the interpolatory projector-splitting integrator from time $t_0$ to $t_1 = t_0 + \Delta t$ can be implemented by sweeping through the cores of $Y$ updating individual cores from $t_0$ to $t_1$. As we update the TT-cores we also update the multi-index sets $\{\I^{\leq j},\I^{>j}\}$ to ensure the interpolatory projectors remain well-defined. We begin with an orthogonal TT representation of solution at time $t_0$ of the form 
%
%
%
\begin{equation}
    Y(t_0) = U_1(t_0) S_1(t_0) V_{>1}(t_0), 
\end{equation}
and the multi-indices $\I^{>j}(t_0)$ defining interpolatory projectors onto the bases $V_{>j}(t_0)$ for $j=1,2,\ldots,d-1$. 
The sweeping algorithm solves the equations in \eqref{split_diffs} sequentially by updating the solution TT-cores $U_j(t_0)$ to $U_j(t_1)$ and then computes the indices $\I^{\leq j}(t_1)$ for the oblique projectors \eqref{oblique_range_co} onto the updated bases $\bm U_{\leq j}(t_1)$ required for the next step in the sweep. 

To begin we apply Theorem \ref{thm:one_core_updates} to solve the first differential equation in \eqref{split_diffs} by integrating 
\begin{equation}
    \frac{d K_1(t)}{dt} = 
     G_1^+\left(:,\I^{>1}(t_0),t\right) \left[V_{>1}\left(:,\I^{>1}(t_0),t_0\right)\right]^{-\top}, \qquad K_1(t_0) = U_1(t_0) \bm S_1(t_0),
\end{equation}
from $t_0$ to $t_1$. 
The solution $Y_1^+(t_1) = K_1(t_1) V_{>1}(t_0)$ is the starting value $Y_1^-(t_0) = Y_1^+(t_1)$ for the second equation in \eqref{split_diffs}. 
We then prepare $Y_1^-(t_0)$ for the application of Theorem \ref{thm:one_core_updates} by decomposing 
$K_1(t_1) = U_1(t_1) \bm R_1(t_1)$ to obtain 
$
Y_1^-(t_0) = U_1(t_1) \bm R_1(t_1) V_{>1}(t_0) 
$
where $U_1(t_1)$ is left-orthogonal and compute indices $\I^{\leq 1}(t_1) = \texttt{DEIM}(\bm U_{\leq 1}(t_1))$. 
Now we can apply Theorem \ref{thm:one_core_updates} to solve the second differential equation in \eqref{split_diffs} by integrating 
\begin{equation}
    \frac{d \bm S_1(t)}{dt} 
    = -\left[U_{\leq 1}\left(\I^{\leq 1}(t_1),:,t_1\right)\right]^{-1} 
G_1^-\left(\I^{\leq 1}(t_1),\I^{>1}(t_0),t\right) 
\left[ V_{>1}\left(:,\I^{>1}(t_0),t_0\right)\right]^{-\top}, \quad \bm S_1(t_0) = \bm R_1(t_1), 
\end{equation}
from time $t_0$ to $t_1$ to obtain the solution 
$Y_1^-(t_1) = U_{1}(t_1) \bm S_1(t_1) V_{>1}(t_0)$. 
The algorithm proceeds recursively with step $j$ of the sweep described below. 

\paragraph{Computation of $Y_j^+(t_1)$} 

The starting value $Y_j^+(t_0) = Y_{j-1}^-(t_1)$ is available in the form 
\begin{equation}
\label{Y1_plus0}
Y_j^+(t_0) = U_{\leq j-1}(t_1) \bm S_{j-1}(t_1) V_{>j-1}(t_0),
\end{equation}
from the computation of $Y_{j-1}^-(t_1)$, 
as are the multi-index sets $\I^{\leq j-1}(t_1)$ and $\I^{>j}(t_0)$. 
To apply Theorem \ref{thm:one_core_updates} we write \eqref{Y1_plus0} as 
\begin{equation}
    Y_j^+(t_0) = U_{\leq j-1}(t_1) \left[\bm S_{j-1}(t_1) V_j(t_0)\right] V_{>j}(t_0), 
\end{equation}
and then integrate 
\begin{equation}
\begin{split}
    \frac{d K_{j}(t)}{dt} = 
    \left[U_{\leq j-1}\left(\I^{\leq j-1}(t_1),:,t_1\right)\right]^{-1} 
         G_j^+\left(\I^{\leq j-1}(t_1),:,\I^{>j}(t_0),t\right) 
        \left[V_{>j}\left(:,\I^{>j}(t_0),t_0\right)\right]^{-\top},\\
        \qquad K_j(t_0) = \bm S_{j-1}(t_1) V_j(t_0),
\end{split}
\end{equation}
from $t_0$ to $t_1$ to obtain the solution 
$Y_j^+(t_1) = U_{\leq j-1}(t_1) K_j(t_1) V_{>j}(t_0)$. 

\paragraph{Computation of $Y_j^-(t_1)$} 
The starting value $Y_j^-(t_0) = Y_j^+(t_1)$ is available 
in the form 
\begin{equation}
\label{Yjminus0}
Y_j^-(t_0) 
= U_{\leq j-1}(t_1) K_j(t_1) V_{>j}(t_0), 
\end{equation}
from the computation of $Y_j^+(t_1)$ as are the multi-index sets $\I^{\leq j-1}(t_1)$ and $\I^{>j}(t_0)$. 
We prepare $Y_j^-(t_0)$ for the application of Theorem \ref{thm:one_core_updates} by decomposing
$K_j(t_1) = U_j(t_1) \bm R_j(t_1)$ 
where $U_j(t_1)$ is left-orthogonal (see Section \ref{sec:orth_TT}) which allows us to write the starting value \eqref{Yjminus0} as 
\begin{equation}
\begin{split}
    Y_j^-(t_0) 
    = U_{\leq j}(t_1) \bm R_j(t_1) V_{>j}(t_0). 
\end{split}
\end{equation}
Then we obtain the multi-indices 
$\I^{\leq j}(t_1)$ from $\I^{\leq j-1}(t_1)$ and the TT-cores $U_{\leq j}(t_1)$ with a substep of the TT-cross-DEIM algorithm as described in \eqref{Wk}-\eqref{Ik}. 
Then by Theorem \ref{thm:one_core_updates}, integrating 
\begin{equation}
    \frac{d \bm S_j(t)}{dt} 
        = -\left[U_{\leq j}\left(\I^{\leq j}(t_1),:,t_1\right)\right]^{-1} 
        G_j^-\left(\I^{\leq j}(t_1),\I^{>j}(t_0),t\right) \left[V_{>j}\left(:,\I^{>j}(t_0),t_0\right)\right]^{-\top}, \quad \bm S_j(t_0) = \bm R_j(t_1), 
\end{equation}
from time $t_0$ to $t_1$ yields the solution 
$Y_j^-(t_1) = U_{\leq j}(t_1) S_j(t_1) V_{>j}(t_0)$. 

Iterating these steps until we obtain $Y_d^+(t_1) = U_{\leq d-1}(t_1) K_d(t_1) = Y(t_1)$ completes one step of the first-order splitting integrator. 
To take another time step the TT representation of $Y(t_1)$ must be orthogonalized from right to left to obtain an orthogonal representation of the solution at time $t$ in the form of \eqref{orth_TT} with $k=1$. 
During this orthogonalization procedure, the indices $\I^{> j}(t_1)$ can be computed with the right-to-left TT-cross-DEIM sweep as described in Section \ref{sec:tt-cross-deim}. 
Similar to the orthogonal projector-splitting integrator, obtaining the second-order Strang projector-splitting integrator is straightforward by composing the Lie-Trotter integrator with its adjoint. 
In this case the forward sweep described above is performed with step-size $\Delta t/2$ and is then followed by a backward sweep also with step-size $\Delta t/2$. 
The oblique projector-splitting integrator has the same computational complexity as the TT-cross integrator. Just as with the corresponding matrix integrators described in Section \ref{sec:matrix_DLR}, the difference between these two integrators is the order in which interpolatory projection and time integration are performed.

\subsection{Rank-adaptive time integration}
\label{sec:rank-adaptive}

The solution to \eqref{nonlinear_ODE0} is often not accurately represented on a tensor manifold $\M_{\r}$ with constant rank for all $t \in [0,T]$. 
Therefore the dynamical low-rank integrators must be able to decrease or increase the solution rank during time integration. 
To decrease the solution rank we use the TT-SVD truncation algorithm at each time $t$ which requires $d-1$ orthogonal representations \eqref{orth_TT} of the solution. 
Such orthogonalizations are required for the TT-cross-DEIM index selection algorithm and thus rank decrease can be performed during time integration with either the TT-cross or interpolatory projector-splitting algorithms at no additional computational cost.

To increase the $k$th component of the TT solution rank during integration with the TT-cross integrator we modify Algorithm \ref{alg:TT_DEIM} to sample $\hat{r}_k>r_k$ indices $\bm l_{\leq k}$ from the left singular vectors \eqref{lk} and $\hat{r}_k$ indices $\bm l_{>k}$ from right singular vectors in \eqref{lk_2} by augmenting the DEIM indices with additional indices selected by another sparse index selection algorithm, e.g., GappyPOD+E \cite{gappy_podE}. 
From the $\bm l_{\leq k}, \bm l_{>k}$ we construct $\I^{\leq k},\I^{>k}$ in \eqref{Ik},\eqref{I_g_k} each with $\hat{r}_k$ indices. 
We then integrate the solution $Y(t)$ forward in time on the manifold $\M_{\hat{\bm r}}$ using the equations \eqref{propagator_interp}. 
It is well-known that the solution $Y(t)$ with rank $\bm r$ belongs to the boundary of the higher rank manifold $\M_{\hat{\bm r}}$ where the tangent space is not well-defined \cite{h_tucker_geom}. 
Nevertheless, the evolution equations \eqref{propagator_interp}, which define the interpolatory tangent space projection, are well-defined on the boundary of $\M_{\hat{\bm r}}$. 
These equations allow us to integrate $Y(t)$ forward in time on $\M_{\hat{\bm r}}$ thereby increasing the solution rank. 
To increase the $k$th component of the TT solution rank during integration with the projector-splitting integrator we add new (orthogonal) basis vectors to the TT cores with zero singular value and then sample indices from this augmented basis and apply the projector-splitting integrator to the augmented solution. 
Once again adding new basis vectors with zero singular values places the approximate solution on the boundary of a higher rank manifol $\M_{\hat{\bm r}}$. 
The projector-splitting integrator is robust to zero singular values and allows us to integrate the solution off of the boundary of the low-rank manifold. 
A simple criterion for determining when to increase the $k$th component of the TT-rank vector is based on the singular values $\{\sigma_k(\alpha_k,t)\}_{\alpha_k=1}^{r_k}$ of the unfolding matrix $\bm Y^{\langle k \rangle}$. 
We select the rank to ensure that the relative size of the smallest singular value 
\begin{equation}
\label{TT_rank_proxy}
\epsilon_k(t) = \frac{\sigma_k(r_k,t)}{ \sqrt{\displaystyle \sum_{\alpha_k=1}^{r_k} \sigma_k(\alpha_k,t)^2}}, \qquad k=1,2,\ldots,d-1 
\end{equation}
remains in a desired range $\epsilon_l \leq \epsilon_k(t) \leq \epsilon_u$. 
This criterion is an adaptation of the rank-adaptive criterion proposed in \cite{Oblique-CUR} for matrix differential equations and subsequently generalized to the Tucker format \cite{ghahremani2024deim}, to the TT format. 

\section{Numerical examples} 
\label{sec:numerics} 

We now apply the proposed dynamical low-rank collocation methods to several tensor differential equations \eqref{nonlinear_ODE0} arising from the discretization of partial differential equations \eqref{nonlinear_PDE0} and compare the accuracy and efficiency with existing time integration schemes on tensor manifolds. 
We measure the accuracy of the low-rank approximations $Y(t)$ to the solution $X(t)$ of \eqref{nonlinear_ODE0} in the relative Frobenius norm 
\begin{equation}
    E(t) = \frac{\left\| Y(t) - X(t) \right\|_F}{\left\| X(t) \right\|_F}. 
\end{equation}
We compute a reference solution $X(t)$ for each application below by integrating the differential equation \eqref{nonlinear_ODE0} with the four-stage explicit Runge-Kutta (RK4) method using time step-size $\Delta t = 10^{-3}$. 

\subsection{2D Vlasov-Poisson equation}

We begin with a two-dimensional example ($d=2$) demonstrating the proposed methods on low-rank matrix manifolds described in Section \ref{sec:matrix_DLR}. We consider the Vlasov-Poisson equation 
\begin{equation}
    \label{VP}
    \begin{cases}
    \displaystyle\frac{\partial u(x,v,t)}{\partial t} + v \displaystyle\frac{\partial u(x,v,t)}{\partial x} + E(x) \displaystyle\frac{\partial u(x,v,t)}{\partial v} = 0 \\
    u(x,v,0) = u_0(x,v),
    \end{cases}
\end{equation}
from \cite[Example 4.4]{guo2022low} 
with initial condition $f(x,v,0) = \exp(-20(x^2 + v^2))$, electric field $E(x) = 0.5\sin(\pi x)$ and $x \in \Omega_x = [-1,1], v \in \Omega_v = [-1,1]$. 
Discretizing $\Omega_x$ and $\Omega_v$ using $n=64$ points and approximating derivatives with a Fourier pseudo-spectral method \cite{spectral_methods_book} we obtain a semi-discrete version of the Vlasov-Poisson equation \eqref{VP} in the form of a differential equation \eqref{nonlinear_ODE0} with $d=2$, i.e., a matrix differential equation. 

\begin{figure}[!t]
\centerline{\footnotesize\hspace{1.5cm} (a)   \hspace{7.5cm} (b) \hspace{1cm}}
\centering
\includegraphics[scale=0.41]{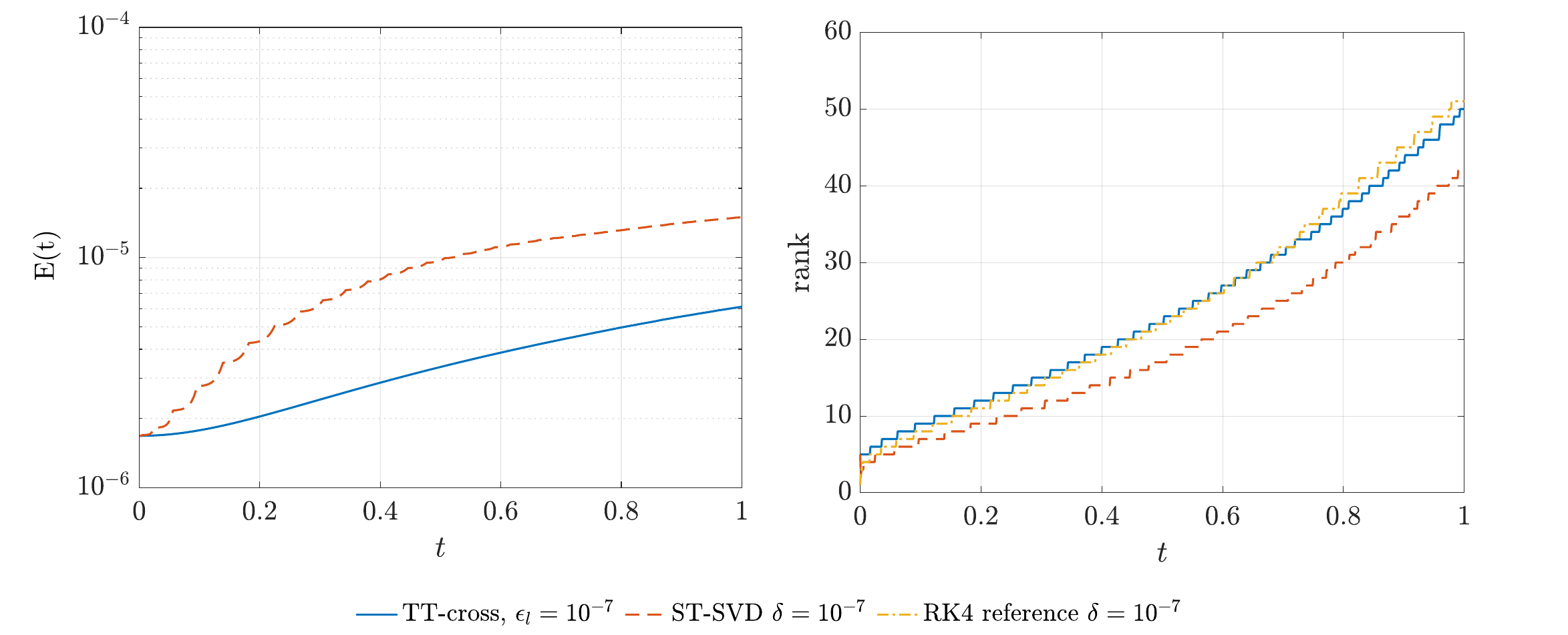}
\centerline{\footnotesize\hspace{2cm} (c)   \hspace{7.5cm} (d) \hspace{1.4cm}}
\includegraphics[scale=0.41]{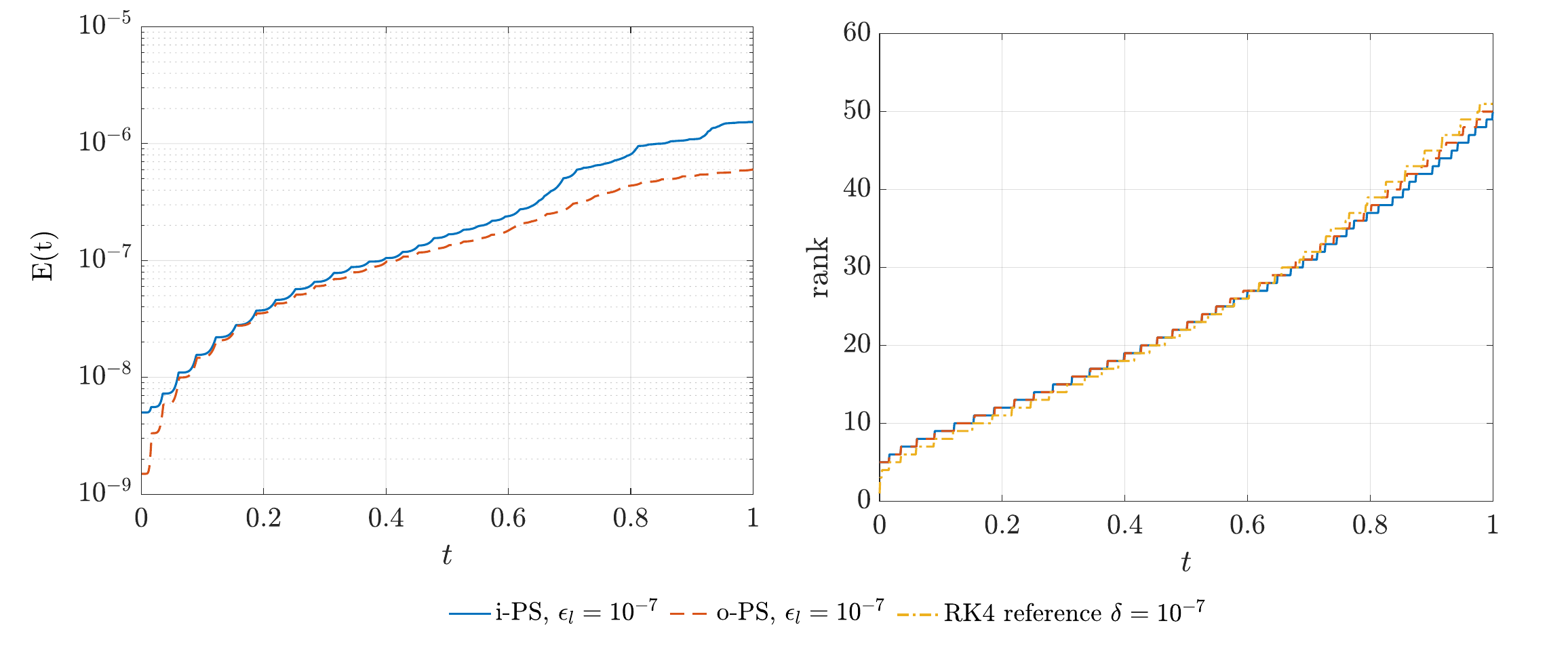}
\caption{
Low-rank approximations to the solution of the two-dimensional Vlasov-Poisson equation \eqref{VP}. 
(a) Relative error versus time of TT-cross and ST-SVD solutions. The TT-cross solution was computed using rank-adaptive singular value threshold $\epsilon_l = 10^{-7}$ and the ST-SVD solution was computed using truncation threshold $\delta = 10^{-7}$. 
(b) Rank versus time of the rank-adaptive TT-cross and ST-SVD solutions and the numerical rank of the reference RK4 solution with singular value threshold $10^{-7}$. 
(c) Relative error versus time of solutions computed with interpolatory and orthogonal projector-splitting using rank-adaptive singular value threshold $\epsilon_l = 10^{-7}$. 
(d) Rank versus time of the rank-adaptive solutions computed with interpolatory and orthogonal projector splitting integrators and the numerical rank of the reference RK4 solution with singular value threshold $10^{-7}$.
}
\label{fig:VP2D}
\end{figure}

We compared the TT-cross integrator presented in Section \ref{sec:TT_cross_int} with a step-truncation method using SVD-based truncation (ST-SVD). 
For both integrators we used Adams-Bashforth 2 with step-size $\Delta t = 10^{-3}$. 
We utilized the rank-adaptive mechanism described in Section \ref{sec:rank-adaptive} for TT-cross with parameter $\epsilon_l = 10^{-7}$. 
For the ST-SVD solution we set relative truncation tolerance $\delta = 10^{-7}$ at each time step allowing the solution rank to adapt in time accordingly. 
In Figure \ref{fig:VP2D}(b) we plot the rank of the TT-cross and ST-SVD solutions versus time and the numerical rank of the reference RK4 solution with singular value threshold $\delta=10^{-7}$, i.e., the number of singular values with relative size larger than $\delta$. 
The rank grows rapidly during time integration which allows us to assess the robustness of the rank-adaptive mechanism for the TT-cross integrator. 
In Figure \ref{fig:VP2D}(a) we plot the relative error of the TT-cross and ST-SVD solutions in the Frobenius norm versus time.
We observe that the TT-cross solution is more accurate than the ST-SVD solution due to the TT-cross solution rank being slightly larger than the rank of the ST-SVD solution at each step. 
The error of the TT-cross solution remains controlled during time integration, demonstrating the effectiveness of the rank-adaptive mechanism.

Next we compared the interpolatory projector-splitting integrator (i-PS) presented in Section \ref{sec:proj_split} with the orthogonal projector-splitting integrator (o-PS) introduced in \cite{lubich2014projector}. 
For the interpolatory and orthogonal projector-splitting integrators we used step-size $\Delta t = 10^{-3}$ and solved the differential equations in the \textbf{K-}, \textbf{S-}, and \textbf{L-step} with RK4. 
We also used the rank-adaptive mechanism described in Section \ref{sec:rank-adaptive} with parameter $\epsilon_l = 10^{-7}$ for both solutions. 
In Figure \ref{fig:VP2D}(b) we plot the solution ranks versus time and the numerical rank of the reference RK4 solution with singular value threshold $10^{-7}$. 
Both solutions have the same rank until approximately $t=0.7$ when the i-PS solution rank becomes slightly smaller than the o-PS solution rank. 
In Figure \ref{fig:VP2D}(a) we plot the relative errors in the Frobenius norm versus time. 
The error of the i-PS and o-PS solutions is similar until around $t=0.4$, at which point the i-PS solution becomes slightly less accurate. 
This difference in accuracy is due to the i-PS method computing a quasi-optimal projection onto the tangent space, while the o-PS method computes the optimal projection at each time step. 
In addition, the slight difference in rank of the solutions also contributes to the difference in accuracy.

\subsection{3D Allen-Cahn equation}
\begin{table}[t]
\caption{CPU-time and accuracy of low-rank methods for integrating the 3D Allen-Cahn equation \eqref{Allen-Cahn}. The ranks were chosen using $\delta =10^{-3}$ and $\delta=10^{-4}$.} 
\label{table:timings}
\begin{center}
\begin{tabular}{ l | c | c | c }
\hline
\textbf{Method} & \textbf{Average rank} $\|\bm r(t)\|_1$ & \textbf{Runtime (min)} & \textbf{Relative Error ($t=10$)} \\ 
\hline
{TT-cross AB2} & 24.2 & 4.0 & $2.5 \times 10^{-2}$ \\
{ST-SVD AB2} & 24.2 & 16.7 & $7.13 \times 10^{-3}$ \\ 
{i-PS RK4} & 24.2 & 23.9 & $2.21 \times 10^{-2}$ \\
{o-PS RK4} & 24.2 & 287.6 & $7.13 \times 10^{-3}$ \\
\hline
{TT-cross AB2} & 32.4 & 4.3 & $3.6 \times 10^{-3}$ \\
{ST-SVD AB2} & 32.4 & 27.2 & $1.0 \times 10^{-3}$ \\
{i-PS RK4} & 32.4 & 21.8 & $3.6 \times 10^{-3}$ \\
{o-PS RK4} & 32.4  & 522.1 & $1.0 \times 10^{-3}$\\
\hline
\end{tabular}
\end{center}
\end{table}

The Allen-Cahn equation is a reaction-diffusion PDE that models phase separation in multi-component alloy systems \cite{Allen_Cahn,Trefethen_2005}. 
A simple form of such equation features a Laplacian and a cubic non-linearity 
\begin{equation}
\label{Allen-Cahn}
\begin{cases}\displaystyle
\frac{\partial u(\bm x,t)}{\partial t} = \alpha \Delta u(\bm x,t) + u(\bm x,t) - u(\bm x,t)^3, \\
u(\bm x,0) = u_0(\bm x). 
\end{cases}
\end{equation}
We consider the spatial domain $\Omega = [0,2\pi]^3$ with periodic boundary conditions, initial condition $u_0(x_1,x_2,x_3) = g(x_1,x_2,x_3) - g(2x_1,x_2,x_3) + g(x_1,2x_2,x_3) - g(x_1,x_2,2x_3)$
where
\begin{equation}
g(x_1,x_2,x_3) = 
\frac{\left(e^{-\tan(x_1)^2} + e^{-\tan(x_2)^2} + e^{-\tan(x_3)^2}\right)
\sin(x_1 + x_2 + x_3) }
{1 + e^{|\csc(-x_1/2)|} + e^{|\csc(-x_2/2)|} + e^{|\csc(-x_3/2)|}}, 
\end{equation}
and diffusion parameter $\alpha = 0.1$. 
Discretizing $\Omega$ using $n=64$ points in each dimension and approximating derivatives with a Fourier pseudo-spectral method \cite{spectral_methods_book}, we obtain a semi-discrete version of the Allen-Cahn equation in the form of \eqref{nonlinear_ODE0}. 
%
%

%
\begin{figure}[!t]
\centerline{\footnotesize\hspace{1.5cm} (a)   \hspace{7.8cm} (b) \hspace{1cm}}
\centering
\includegraphics[scale=0.41]{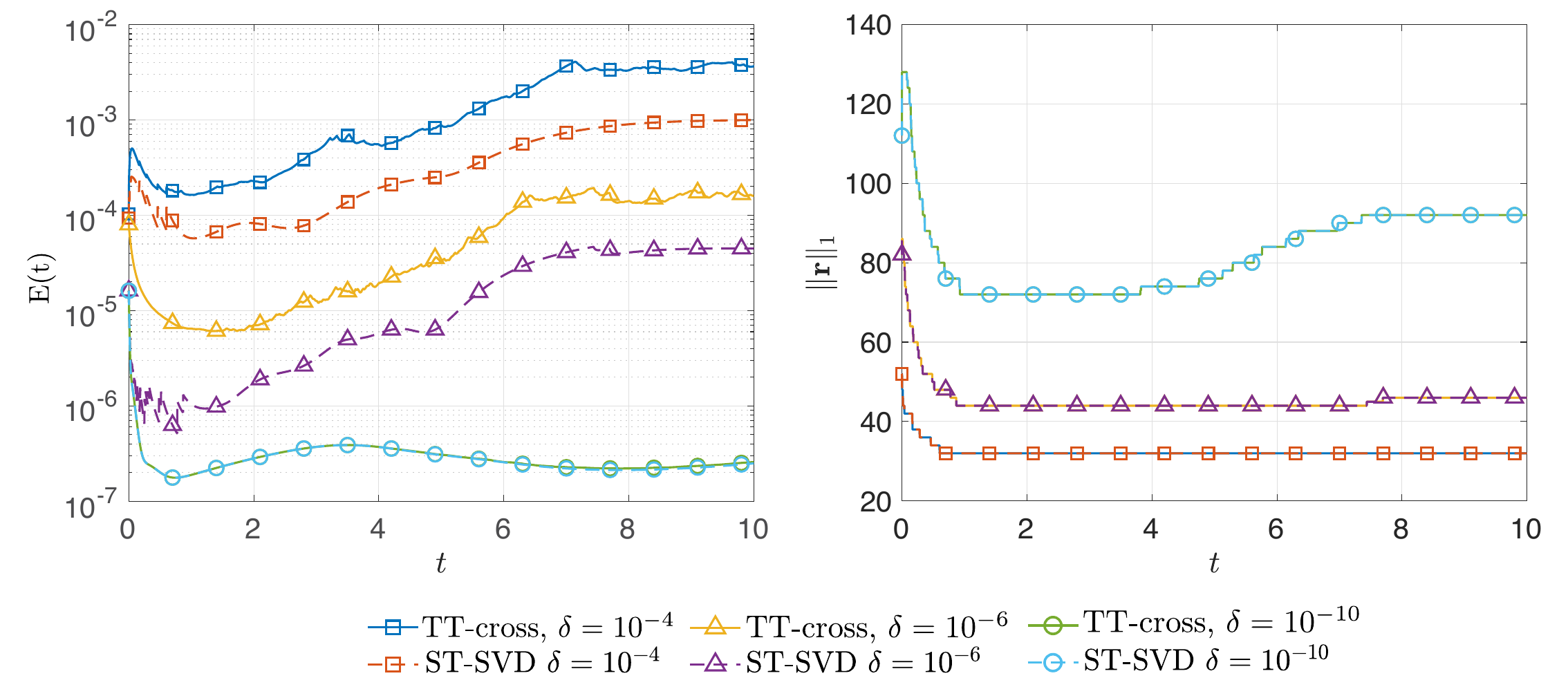}
\centerline{\footnotesize\hspace{2cm} (c)   \hspace{8cm} (d) \hspace{1cm}}
\includegraphics[scale=0.41]{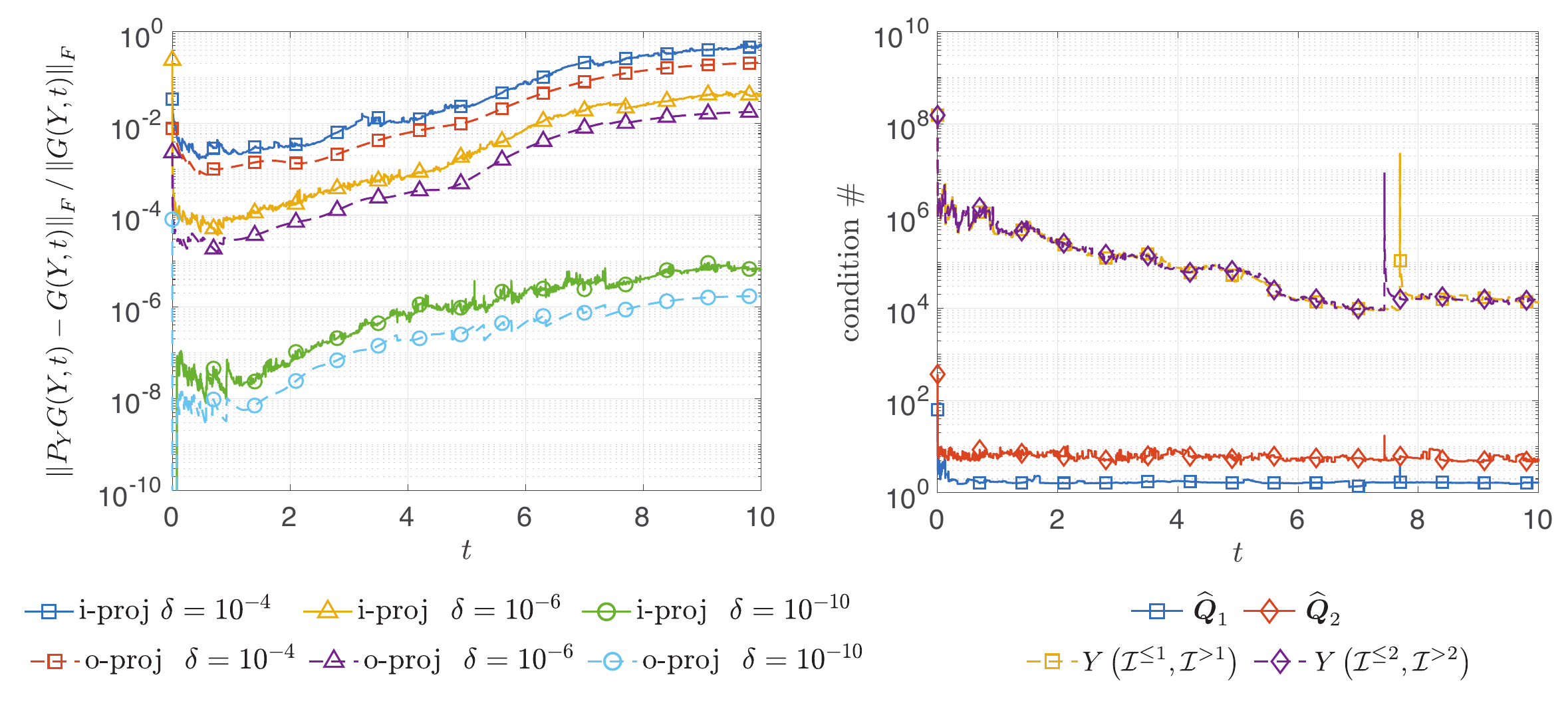}
\caption{
Low-rank approximations to the solution of the three-dimensional Allen-Cahn equation \eqref{Allen-Cahn} computed with the TT-cross and ST-SVD methods. 
The ranks were determined using different truncation tolerances $\delta = 10^{-4},10^{-6},10^{-10}$ in the ST-SVD method.  
(a) Relative error in the Frobenius norm versus time. 
(b) 1-norm of the TT-rank vector versus time. 
(c) Relative error of interpolatory (i-proj) and orthogonal (o-proj) projections onto the tensor manifold tangent space versus time. 
(d) Condition number of the matrices non-orthogonalized matrices in \eqref{TT_cross_time} and the corresponding orthogonalized matrices in \eqref{TT_cross_time2} used to construct the TT-cross solution at each time step. 
}
\label{fig:AC3D}
\end{figure}

We compared the TT-cross integrator presented in Section \ref{sec:TT_cross_int} with the step-truncation SVD (ST-SVD) integrator \cite{rodgers2020step-truncation} using different relative truncation tolerances $\delta = 10^{-3},10^{-4},10^{-6},10^{-10}$ for determining the solution rank at each time step. 
We set the solution rank in the TT-cross simulations equal to the ranks obtained from the ST-SVD simulations with truncation tolerances in order to compare the methods for solutions computed with the same rank. 
The rank decrease was performed using TT-SVD truncation and the rank increase by sampling more tensor cross indices than singular vectors using the GappyPOD+E algorithm \cite{gappy_podE} as described in Section \ref{sec:rank-adaptive}. 
Time integration for both ST-SVD and TT-cross was performed with Adams-Bashforth 2 and step-size $\Delta t = 10^{-3}$. 

In Figure \ref{fig:AC3D}(b), we plot the 1-norm of the ST-SVD and TT-cross solution ranks. 
The smoothing effects due to diffusion in the Allen-Cahn equation cause the TT ranks to decay relatively quickly from time $t=0$ to time $t\approx 1.5$. 
In Figure \ref{fig:AC3D}(a), we plot the relative error measured in the Frobenius norm of the ST-SVD and TT-cross solutions versus time. 
The ST-SVD solution is more accurate than the TT-cross solution computed with the same rank, which is expected. 
Indeed, the ST-SVD method computes the best rank-$\bm r$ projection of the solution onto the low-rank manifold $\M_{\r}$ at each time step while the TT-cross method computes a quasi-optimal projection onto the tangent space of the manifold at each time step. 
When the rank of the TT solutions is large enough (in this case corresponding to $\delta = 10^{-10}$), the time integration error dominates the low-rank approximation error and the ST-SVD and TT-cross methods produce solutions with the same accuracy. 
When the low-rank error dominates the time integration error ($\delta = 10^{-4},10^{-6}$) we observe in Figure \ref{fig:AC3D}(a) that the ST-SVD is about half an order of magnitude more accurate than the TT-cross solution of the same rank for all ranks and at each time $t$.  
In Figure \ref{fig:AC3D}(c), we compare the accuracy of the interpolatory projection (i-proj) onto the tangent space computed from the TT-cross solution and the orthogonal projection (o-proj) onto the tangent space computed from the ST-SVD simulation. 
The orthogonal projection is more accurate than the interpolatory projection by approximately one order of magnitude or less at each time $t$. 
Similar to the difference in error between the solutions, the difference in error between the i-proj and o-proj is constant over all ranks and for all time $t$. 

\begin{figure}[!t]
\centerline{\footnotesize\hspace{1.8cm} (a)   \hspace{7.8cm} (b) \hspace{1cm}}
\centering
\includegraphics[scale=0.41]{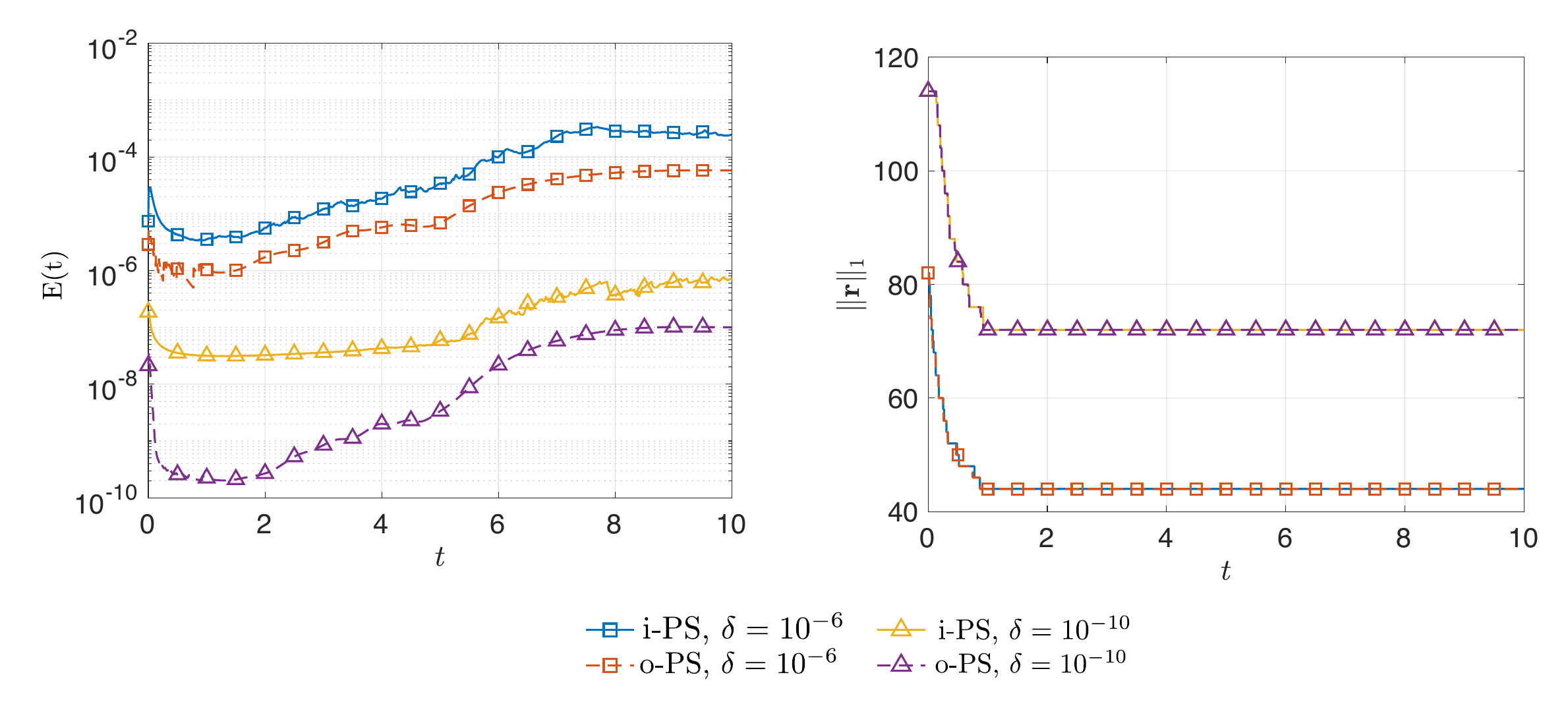}
\caption{
Low-rank approximations to the solution of the three-dimensional Allen-Cahn equation \eqref{Allen-Cahn} computed with the interpolatory and orthogonal projector-splitting integrators. 
Solutions are truncated at each time $t$ using tolerances TT-SVD with relative tolerance $\delta = 10^{-6},10^{-10}$. 
(a) Relative error in the Frobenius norm versus time. 
(b) 1-norm of the TT-rank vector versus time. 
}
\label{fig:AC3D_split}
\end{figure}

The improved accuracy of the ST-SVD method over the TT-cross method comes at a significant computational cost due to the cubic nonlinearity in the Allen-Cahn equation \eqref{Allen-Cahn}. 
The reason is that the ST-SVD method requires computing a TT representation of $G(Y(t),t)$ at each time $t$, which is costly. 
Indeed, recall that standard algorithms for multiplying two TTs $Y_1$ and $Y_2$ with ranks $\bm r_1=\begin{bmatrix} r_1 & \cdots & r_1 \end{bmatrix}$ and $\bm r_2 = \begin{bmatrix} r_2 & \cdots & r_2\end{bmatrix}$ results in a TT $Y_1 Y_2$ with rank equal to the Hadamard (element-wise) product of the two ranks $\bm r_1 \circ \bm r_2$. 
These ranks are in general not optimal and to control the TT rank we perform a TT-SVD truncation requiring $\mathcal{O}(dn (r_1 r_2)^3)$ operations. 
We used two TT-SVD truncations $\mathfrak{T}_{\delta}$ with relative accuracy $\delta$ to compute the cubic term 
\begin{equation}
\label{cubic_truncation}
(Y)^3 = \mathfrak{T}_{\delta}^{\rm svd}\left( Y \  \mathfrak{T}_{\delta}^{\rm svd}\left(Y Y \right) \right),
\end{equation}
incurring a cost of $\mathcal{O}(dnr^6)$ operations at each time $t$. 
It is possible to accelerate the computation $G(Y,t)$ by carrying out sums and products of TTs with approximate low-rank tensor arithmetic, black-box tensor cross approximation \cite{parr_TT_cross_Dolgov}, or randomized algorithms \cite{Daas_randomized_TT_rounding}. 
However such algorithms introduce additional errors in the low-rank approximation that can be difficult to control. 
In comparison, the TT-cross integrator does not require $G(Y,t)$ in a low-rank form and instead evaluates $G(Y,t)$ at $\mathcal{O}(d n r^2)$ indices. 
Thus the computational cost of the cubic nonlinearity for TT-cross is negligible compared to the $\mathcal{O}(dnr^3)$ cost of the TT-cross-DEIM index selection algorithm and evaluating the subtensors of $Y(t)$ required to integrate the system of equations \eqref{propagator_interp}. 
%

We also compared the interpolatory projector-splitting (i-PS) integrator presented in Section \ref{sec:proj_split} with the orthogonal projector-splitting (o-PS) integrator from \cite{Lubich_2015} using two different truncation tolerances $\delta = 10^{-6},10^{-10}$ on the singular values of the solutions. 
In both cases we used first-order Lie-Trotter splitting with time step-size $\Delta t = 10^{-3}$ and solved each of the substeps in \eqref{split_diffs} with RK4. 
In Figure \ref{fig:AC3D_split}(a) we plot the error of the solutions computed with the i-PS and o-PS methods versus time. 
We observe that the i-PS method is less accurate than the o-PS method. This is expected since the i-PS method integrates $Y(t)$ on $\M_{\r}$ using a quasi-optimal tensor in the tangent space while the o-PS method integrates uses the optimal tensor in the tangent space. 
The difference in error is similar to the comparison of TT-cross and ST-SVD except for $t \in [0,5]$ in the simulations using $\delta=10^{-10}$ where the difference in error is significantly larger. 
In Figure \ref{fig:AC3D_split}(b) we plot the ranks of the i-PS and o-PS solutions versus time. 

In Table~\ref{table:timings} we compare the runtime and relative error at time $t=10$ of the low-rank solutions computed with existing methods (ST-SVD AB2, o-PS RK4) with the solutions computed using the proposed methods (TT-cross AB2, i-PS RK4). We consider two different rank-adaptive simulations with ranks determined by $\delta =10^{-3},10^{-4}$ and report the average $1$-norm of the rank vector over all time steps. 
The TT-cross AB2 method with an average rank of $24.2$ is approximately $4.2$ times faster than the ST-SVD AB2 method at the same rank, while being roughly half an order of magnitude less accurate. With an average rank of $32.4$, the TT-cross AB2 method is approximately $6.3$ times faster than the ST-SVD AB2 method, while being less than half an order of magnitude less accurate. 
The speedup observed for the projected RK4 methods is even greater, as these methods require more evaluations of the right-hand side, which includes the cubic nonlinearity. The interpolatory RK4 method with an average rank of $24.2$ is approximately $12$ times faster than the o-PS RK4 method at the same rank, while being roughly half an order of magnitude less accurate. With an average rank of $32.4$, the i-PS RK4 method is approximately $24$ times faster than the ST-SVD AB2 method, while being less than half an order of magnitude less accurate.

\subsection{4D advection-diffusion-reaction equation}

Finally we consider the advection-diffusion-reaction (ADR) equation 
\begin{equation}
\label{ADR}
\begin{cases}
\displaystyle\frac{\partial u(\bm x,t)}{\partial t} = 
\nabla \cdot \left(\bm \mu_i(\bm x,t) u(\bm x,t)\right) 
+ \sigma \Delta u(\bm x,t) + R(u) \\
u(\bm x,0) = u_0(\bm x),
\end{cases}
\end{equation}
where $R(u)$ is a nonlinear reaction term. 
We consider the spatial domain $\Omega = [0,2\pi]^4$ with periodic boundary conditions and set 
\begin{equation}
    p_0(\bm x) = 
    \exp(\sin(x_1)\sin(x_2)\sin(x_3)\sin(x_4)),
\end{equation}
$R(u) = -0.1 u/(1+u^2)$, $\sigma = 1/4$ and 
\begin{equation}
\label{drift_diffusion}
\bm \mu(\bm x) = \frac{1}{2} \begin{bmatrix}
g(x_2,x_3) \\
g(x_3,x_4) \\
g(x_4,x_1) \\
g(x_2,x_3) 
\end{bmatrix},
\end{equation}
where $g(x,y)=\exp(\sin(x)\cos(y))$. 
Discretizing $\Omega$ using $n=32$ points in each dimension and approximating derivatives with a Fourier pseudo-spectral method \cite{spectral_methods_book} we obtain a semi-discrete version of the ADR equation \eqref{ADR} in the form of \eqref{nonlinear_ODE0}. 
%
%
%
%
\begin{figure}[!t]
\centerline{\footnotesize\hspace{1.0cm} $t=0$   \hspace{6.85cm} $t=1$ \hspace{1cm}}
\centering
\includegraphics[scale=0.35]{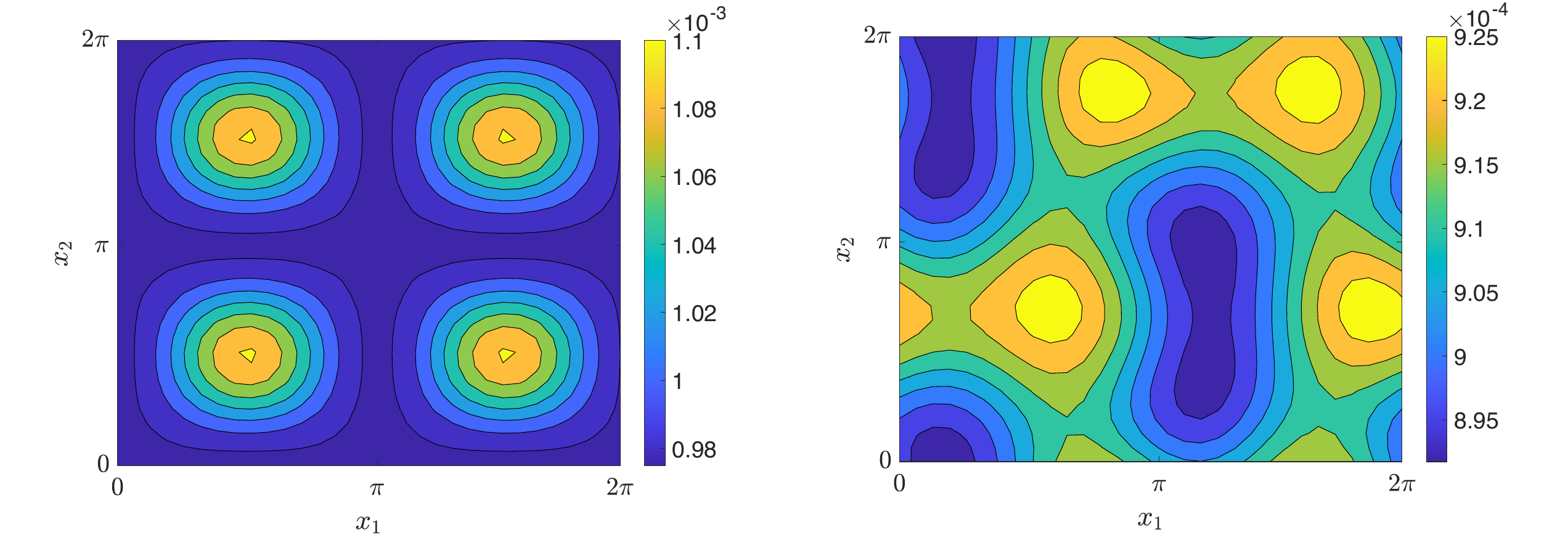}
\caption{
The $(x_1,x_2)$-marginals of the reference solution to the four-dimensional ADR equation \eqref{ADR} at time $t=0$ and $t=1$. 
}
\label{fig:FP4D_marginals}
\end{figure}

\begin{figure}[!t]
\centerline{\footnotesize\hspace{1.7cm} (a)   \hspace{7.5cm} (b) \hspace{1cm}}
\centering
\includegraphics[scale=0.42]{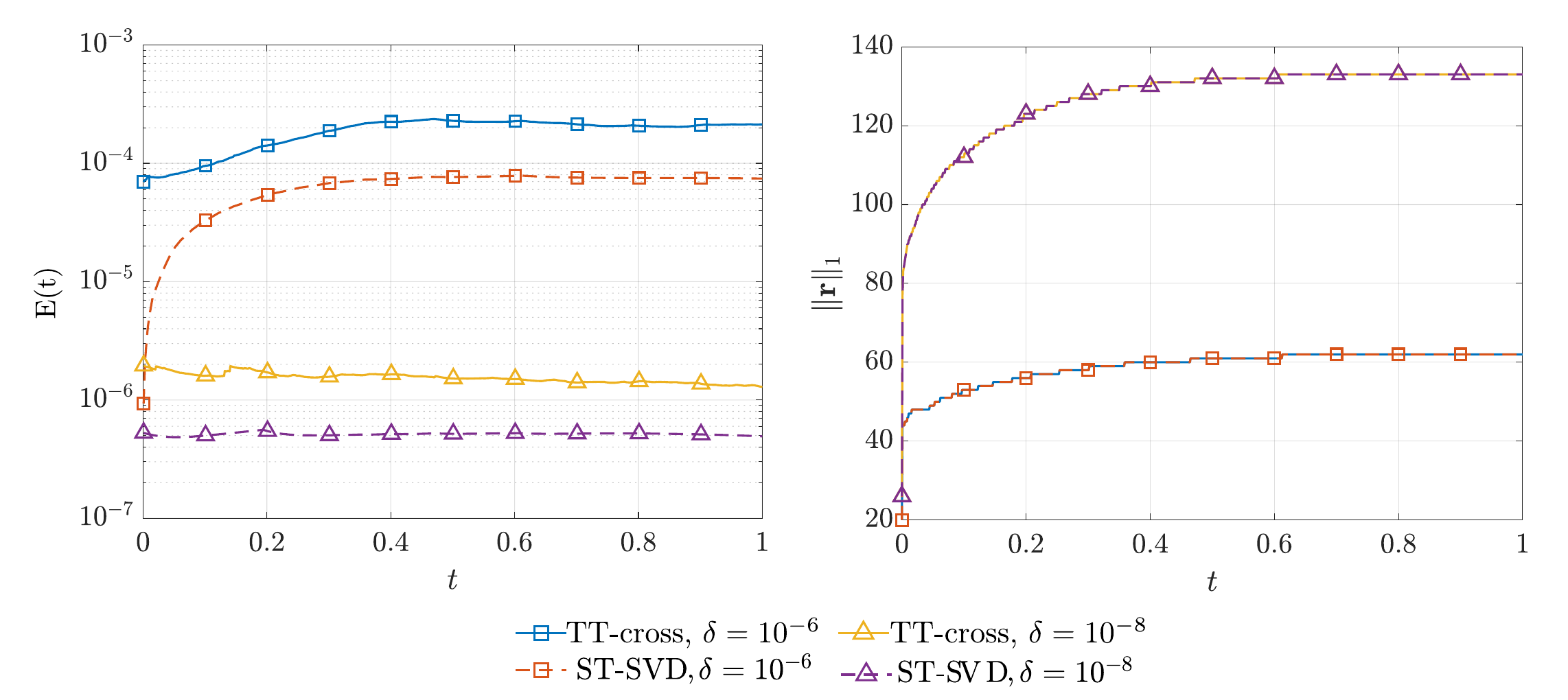}
\caption{
Low-rank approximations to the solution of the four-dimensional ADR equation \eqref{ADR} computed with the TT-cross and ST-SVD integrators. The ranks were determined using different truncation tolerances $\delta = 10^{-6},10^{-8}$ in the ST-SVD method.  
(a) Relative error in the Frobenius norm versus time. 
(b) 1-norm of the TT-rank vector versus time. 
}
\label{fig:AD4_error_rank}
\end{figure}

We computed two approximate low-rank solutions on a TT manifold \eqref{TT-mfld} with the step-truncation SVD method (ST-SVD) \cite{rodgers2020step-truncation} using different relative truncation tolerances $\delta = 10^{-6},10^{-8}$. 
Computing the ST-SVD solution requires a low-rank approximation of $G(Y(t),t)$ at each time $t$, which is challenging for the nonlinear ADR equation \eqref{ADR} as there are no reliable algorithms available for computing the fractional nonlinearity in the low-rank format. 
To compute the $G(Y(t),t)$, we construct the full tensor representation of the TT-SVD solution with $n^4$ degrees of freedom, compute the fractional nonlinearity, and then compress the result into a TT with a recursive SVD. 
This approach is of course not viable in higher dimensions but it allows us to compare our TT-cross solution with the ST-SVD method in this case which computes the best low-rank approximate solution at each time step. 
The map $G$ obtained from discretizing \eqref{ADR} includes four coefficient tensors $c_1,c_2,c_3,c_4\in \mathbb{R}^{n \times n \times n \times n}$ (resulting from the discretization of $g(x,y)$) that are not expressed in a low-rank format upon discretization of $G$. 
In order to compute $G(Y(t),t)$ in low-rank format at each time, we decomposed the four coefficient tensors in $G$ using TT-SVD compression with relative accuracy $\delta$. 
For $\delta = 10^{-6}$ and $\delta = 10^{-8}$ we obtained coefficient tensors of the same rank 
\begin{equation}
\begin{split}
    \text{TT-rank}(c_1) &= \begin{bmatrix}
        1 & 1 & 12 & 1 & 1 
    \end{bmatrix}, \\
        \text{TT-rank}(c_2) &= \begin{bmatrix}
        1 & 1 & 1 & 12 & 1 
    \end{bmatrix}, \\
        \text{TT-rank}(c_3) &= \begin{bmatrix}
        1 & 13 & 13 & 13 & 1 
    \end{bmatrix}, \\
        \text{TT-rank}(c_4) &= \begin{bmatrix}
        1 & 1 & 12 & 1 & 1 
    \end{bmatrix}. \\
\end{split}
\end{equation}
We computed $G(Y(t),t)$ in the ST-SVD method at each time step by taking products of the low-rank approximate coefficient tensors $c_k$ with the low-rank solution tensor $Y$ and then used TT-SVD truncation to compress the product. 
We then added the TT representation of the reaction term and  applied TT-SVD truncation after adding two low-rank tensors in order to control tensor rank when computing $G(Y(t),t)$ at each time $t$. 
Time integration for the ST-SVD simulation was performed with AB2 and time step-size $\Delta t = 10^{-3}$. 
In Figure \ref{fig:AD4_error_rank}(b) we plot the $1$-norm of the TT-rank of each ST-SVD solution versus time. 
We observe that the ranks of both solution increase until around $t=0.5$ and then stabilize for $t \in [0.5,1]$. 

We then computed two approximate low-rank solutions on the TT manifold $\M_{\r}$ using the proposed TT-cross integrator. 
In order to compare the results with the ST-SVD simulations we set the solution ranks in the TT-cross simulations equal to the ranks obtained from the ST-SVD simulations with truncation tolerances $\delta = 10^{-6},10^{-8}$. 
We computed the right-hand side of the TT-cross evolution equations \eqref{propagator_interp} by simply evaluating the coefficient tensors at the indices determined by the TT-cross-DEIM Algorithm at each time step. 
The cost of computing the right hand-hand side for the TT-cross evolution equations is negligible compared to the $\mathcal{O}(dnr^3)$ cost of the TT-cross-DEIM index selection algorithm and evaluating the subtensors of $Y(t)$ required to integrate the system of equations \eqref{propagator_interp}. 
In Figure \ref{fig:AD4_error_rank}(a) we compare the relative error in the Frobenius norm of the TT-cross solutions and the ST-SVD solutions. 
We observe that the TT-cross solutions are less accurate than the ST-SVD solutions of the same rank and the difference in accuracy is constant over all solution ranks and for all time $t$. 
This is expected as the ST-SVD method computes the best rank-$\bm r$ projection of the solution onto the TT manifold $\M_{\r}$ at each time step while the TT-cross method computes a quasi-optimal projection onto the tangent space of the low-rank manifold at each time step. 

\section{Conclusions}
\label{sec:conclusions}

We introduced new general purpose dynamical low-rank methods for solving nonlinear differential equations on low-rank manifolds. 
The methods rely on a particular class of oblique projectors onto the tangent space with a low-rank manifold characterized by a cross-interpolation property. 
Such projectors collocate the differential equation on a low-rank tensor manifold and give rise to efficient time integration schemes that allow us to integrate differential equations defined by vector fields without low-rank structure on low-rank manifolds. 
To construct the oblique projections we introduced a new index selection algorithm based on the DEIM for constructing interpolatory projectors in the TT format. 
Furthermore, we showed that such indices also parameterize low-rank TT manifolds and their tangent spaces with cross interpolation.  
Our numerical results demonstrate that the oblique projections onto the tangent space yield good approximations on the low-rank manifold in the Frobenius norm that are efficiently computed for problems defined by vector fields without low-rank structure. 
Our proposed methods thus make dynamical low-rank approximation applicable to a broader class of differential equations and facilitate its use in various practical applications.

\section*{Declarations of interest} None.  

\section*{Acknowledgements} 
This material is based upon work supported by the U.S. Department of Energy, Office of Science, Office of Advanced Scientific Computing Research, Scientific Discovery through Advanced Computing (SciDAC) program under the contract No. DE-AC02-05CH11231.  This research used resources of the National Energy Research Scientific Computing Center, a DOE Office of Science User Facility supported by the Office of Science of the U.S. Department of Energy under Contract No. DE-AC02-05CH11231 using NERSC award ASCR-ERCAP-m1027.


\newpage
\bibliographystyle{plain}
\bibliography{refs}

\end{document}